\documentclass{article}
\usepackage{amsfonts}
\usepackage{amsmath}

\newtheorem{theorem}{Theorem}

\newtheorem{corollary}[theorem]{Corollary}

\newtheorem{lemma}{Lemma}[section]

\newtheorem{remark}[theorem]{Remark}

\numberwithin{equation}{section}

\begin{document}

\title{Small deviations for critical Galton-Watson processes with infinite
variance }
\author{ Vladimir Vatutin \thanks{%
Steklov Mathematical Institute of Russian Academy of Sciences, 8 Gubkina
St., Moscow 119991 Russia Email: vatutin@mi-ras.ru}, Elena Dyakonova\thanks{%
Steklov Mathematical Institute of Russian Academy of Sciences, 8 Gubkina
St., Moscow 119991 Russia Email: elena@mi-ras.ru}, and Yakubdjan Khusanbaev%
\thanks{%
V.I. Romanovskiy Institute of Mathematics, Uzbekistan Academy of Sciences,
Olmazor district, 9 University St. Tashkent 100174, Uzbekistan Email:
yakubdjank@mail.ru}, }
\date{}
\maketitle

\begin{abstract}
We study the asymptotic behavior of small deviation probabilities for the
critical Galton-Watson processes with infinite variance of the offspring
sizes of particles and apply the obtained result to investigate the
structure of a reduced critical Galton-Watson process.
\end{abstract}


\textbf{MSC-class}: 60J80, 60F10

\footnotetext[0]{%
The work of V.A. Vatutin and E.E. Dyakonova was performed at the Steklov
International Mathematical Center and supported by the Ministry of Science
and Higher Education of the Russian Federation (agreement no.
075-15-2025-303).}

\section{Introduction}

We consider a Galton-Watson branching process $\left\{ Z(n),n\geq 0\right\} $
starting with a single founding ancestor in generation 0. Let $f(s)=\mathbf{E%
}s^{\xi }$ be the probability generating function for the offspring number $%
\xi $ of one particle. We assume that the reproduction law of particles
satisfies the condition
\begin{equation}
\mathbf{E}\xi =1,\qquad f(s)=s+\left( 1-s\right) ^{1+\alpha }L(1-s),\quad
0\leq s\leq 1,  \label{MainAssump}
\end{equation}%
where $\alpha \in (0,1]$ and $L(z)$ is a slowly varying function as $%
z\downarrow 0$. Denote $f_{0}(s)=s$ and introduce iterations of the function
$f(s)$ by the equalities
\begin{equation*}
f_{n}(s)=f(f_{n-1}(s)),n=1,2,\ldots .
\end{equation*}%
It is known (see, for instance, \cite{Sl1968}) that, given (\ref{MainAssump})%
\begin{equation}
(1-f_{n}(0))^{\alpha }L\left( 1-f_{n}(0)\right) \sim \frac{1}{\alpha n}
\label{SimpeForm}
\end{equation}%
as $n\rightarrow \infty $ and, therefore, (see, for instance, \cite[Section
1.5, point $5^{0}$]{Sen76} or \cite{Sen67})
\begin{equation}
Q(n):=\mathbf{P}\left( Z(n)>0\right) =1-f_{n}(0)\sim \frac{L^{\ast }(n)}{%
n^{1/\alpha }},\text{ \ \ \ \ }n\rightarrow \infty ,  \label{SurvivalProbab}
\end{equation}%
where $L^{\ast }(n)$ is a function slowly varying as $x\rightarrow \infty $.
Moreover \cite{Sl1968}, for any $\lambda >0$%
\begin{eqnarray}
\lim_{n\rightarrow \infty }\mathbf{E}\left[ e^{-\lambda
(1-f_{n}(0))Z(n)}|Z(n)>0\right] &=&\int_{0}^{\infty }e^{-\lambda x}dM\left(
x\right)  \notag \\
&=&1-\frac{1}{\left( 1+\lambda ^{-\alpha }\right) ^{1/\alpha }},
\label{Yaglom000}
\end{eqnarray}%
where $M(x)$ is a proper nondegenerate distribution such that $M(x)>0$ for
any $x>0$.

It follows from (\ref{Yaglom000}) that, for any $x>0$
\begin{equation*}
\mathbf{P}\left( 0<(1-f_{n}(0))Z(n)\leq x\right) \approx M(x)\mathbf{P}%
\left( Z(n)>0\right)
\end{equation*}%
as $n\rightarrow \infty $. Thus, the typical size of the population survived
up to a distant moment $n$ is of order $(1-f_{n}(0))^{-1}$.

In this paper we study the asymptotic behavior of small deviation
probabilities for the critical Galton-Watson branching process, i.e., we
investigate the probability
\begin{equation}
\mathcal{H}(n):=\big\{0<(1-f_{\varphi (n)}(0))Z(n)\leq 1\big\}.
\label{Def_H}
\end{equation}%
where $\varphi (n)/n\rightarrow 0$ as $n\rightarrow \infty ,$ and apply the
obtained result to the analysis of the structure of a reduced critical
Galton-Watson process under the condition $\mathcal{H}(n).$

The study of small deviation probabilities for Galton-Watson processes is
closely related to the problem of describing the asymptotic behavior of the
local probabilities $\mathbf{P}(Z(n)=j)$ when $j$ is smaller in order than $%
(1-f_{n}(0))^{-1}$. Probabilities of such a kind for the critical
Galton-Watson processes or continuous time Markov branching processes were
investigated by several authors. Thus, Slack \cite[Theorem 2]{Sl1968} has
shown that if condition~(\ref{MainAssump}) is valid, then the limit
\begin{eqnarray}
&&\lim_{n\rightarrow \infty }L\left( 1-f_{n}(0)\right) ^{1/\alpha }\left(
n\alpha \right) ^{1+1/\alpha }\mathbf{P}\left( Z(n)=j\right)  \notag \\
&&\qquad \qquad =\lim_{n\rightarrow \infty }\frac{n\alpha }{1-f_{n}(0)}%
\mathbf{P}\left( Z(n)=j\right) =\mu _{j}\in \lbrack 0,\infty )
\label{Stattionar}
\end{eqnarray}%
exists for any fixed $j\geq 1$, and
\begin{equation}
\mu _{j}=\sum_{l=1}^{\infty }\mu _{l}\mathbf{P}\left( l,j\right) ,\quad
j\geq 1,\quad \sum_{l=1}^{\infty }\mu _{l}\mathbf{P}^{l}\left( \xi =0\right)
=1,  \label{Stattionar2}
\end{equation}%
where $\mathbf{P}\left( l,j\right) =\mathbf{P}\left( Z(n+1)=j|Z(n)=l\right) $
are the transition probabilities of the branching process. Moreover,
\begin{equation}
\sum_{j=1}^{T}\mu _{j}\sim \frac{1}{\alpha \Gamma \left( 1+\alpha \right) }%
\frac{T^{\alpha }}{L\left( 1/T\right) }  \label{Sl2}
\end{equation}%
as $T\rightarrow \infty $, where $\Gamma (x)$ is the standard Gamma-function.

Asymptotic behavior of local probabilities for the critical processes with
finite variance of the reproduction laws of particles was analyzed in \cite%
{Chist1957} and \cite{Zol57} (for the continuous time Markov branching
processes) and in \cite{KNS66} -- \cite{NW2006} (for the critical
Galton-Watson branching processes). To the authors knowledge, the best
result in this direction is contained in \cite{NW2006} and has the following
form: if%
\begin{equation}
f^{\prime }(1)=1,\quad \sigma ^{2}=f^{\prime \prime }(1)\in (0,\infty ),
\label{FinVar}
\end{equation}%
the distribution of $\xi $ is lattice with span $b$ and parameters $j$ and $%
n $ tend to infinity in such a way that the ratio $jn^{-1}$ remains bounded,
then%
\begin{equation*}
\frac{\sigma ^{4}n^{2}}{4b}\left( 1+\frac{2b}{\sigma ^{2}n}\right) ^{j+1}%
\mathbf{P}\left( Z(n)=bj|Z(0)=1\right) \rightarrow 1.
\end{equation*}%
In particular, if $jn^{-1}\rightarrow 0$, then%
\begin{equation*}
\mathbf{P}\left( Z(n)=bj|Z(0)=1\right) \sim \frac{4b}{\sigma ^{4}n^{2}}\sim
\frac{1-f_{n}(0)}{n}\frac{2b}{\sigma ^{2}},
\end{equation*}%
where the last equivalence is justifies by the asymptotic relation
\begin{equation*}
1-f_{n}(0)\sim \frac{2}{\sigma ^{2}n},\ n\rightarrow \infty ,
\end{equation*}%
valid under the condition (\ref{FinVar}) (see, for instance, \cite{AN72},
Ch.1, Sec.9). Hence it follows that if a parameter $\varphi (n)$ tends to
infinity in such a way that $\varphi (n)/n\rightarrow 0$ as $n\rightarrow
\infty ,$ then
\begin{eqnarray}
\mathbf{P}\left( \mathcal{H}(n)\right) &=&\sum_{0<(1-f_{\varphi
(n)}(0))bj\leq 1}\mathbf{P}\left( Z(n)=bj\right)  \notag \\
&\sim &\sum_{0<j\leq \sigma ^{2}\varphi (n)/\left( 2b\right) }\mathbf{P}%
\left( Z(n)=bj\right) \sim \frac{1-f_{n}(0)}{n}\varphi (n)\text{.}
\label{FinVar22}
\end{eqnarray}

Our first aim is to extend (\ref{FinVar22}) to the case when the variance of
$\xi $ is infinite.

\begin{theorem}
\label{T_deviation} Let condition (\ref{MainAssump}) be valid for $\alpha
\in (0,1)$. If $\varphi (n),\ n=1,2,\ldots ,$ is a deterministic function
such that
\begin{equation}
\varphi (n)\rightarrow \infty \text{ and }\frac{\varphi (n)}{n}\rightarrow 0
\label{Ratio}
\end{equation}%
as $n\rightarrow \infty ,$ then
\begin{equation}
\mathbf{P}\left( \mathcal{H}(n)\right) \sim \frac{1-f_{n}(0)}{\alpha n}\frac{%
1}{\Gamma \left( 1+\alpha \right) }\varphi (n).  \label{LargeDeviation}
\end{equation}
\end{theorem}

We see that the asymptotic relation (\ref{LargeDeviation}) is in full
agreement with (\ref{FinVar22}). Note, however, that (\ref{LargeDeviation})
cannot be deduced from (\ref{Stattionar}), since the last is proved for each
fixed $j$ and not for $j\rightarrow \infty $.

For some reason to be explained later on we do not consider here the case $%
\alpha =1$ and $\sigma ^{2}=\infty $. The reader can find a deteiled
discussion concerning this case at the end of the proof of Theorem \ref%
{T_deviation} (see Remark \ref{R_alpha1}).

Theorem \ref{T_deviation} allows us to study, for each fixed $n$ properties
of the so-called reduced Galton-Watson process $\{Z(m,n),0\leq m\leq n\}$,
where $Z(m,n)$ is the number of particles in the original process $%
\{Z(m),m=0,1,\ldots,n\}$ in generation $m<n$ having a positive number of
descendants in generation $n$ with the natural agreement $Z(n,n)=Z(n)$.

Reduced Galton--Watson branching processes were first considered by
Fleischmann and Prehn \cite{FP}, who discussed the subcritical case. The
distribution of the distance $d(n)$ to the most recent common ancestor
(MRCA) \ of all particles existing in that generation for supercritical
Galton-Watson processes and for the critical Galton-Watson processes with
infinite variance of the number of direct descendants of particles have been
investigated by Zubkov \cite{Zub}. He showed that if condition (\ref%
{MainAssump}) is fulfilled, then

\begin{equation*}
\lim_{n\rightarrow \infty }\mathbf{P}\left( \frac{d(n)}{n}\leq
y|Z(n)>0\right) =y
\end{equation*}%
for all $y\in \lbrack 0,1]$. Sagitov \cite{Sag85} considered the reduced
Galton-Watson branching process $\ \left\{ Z(nt,n),0\leq t\leq 1\right\} $
given $\{Z(n)>0\}$ and proved that if condition (\ref{MainAssump}) is
fulfilled, then, as $n\rightarrow \infty $ the process weakly converges in
Skorokhod topology to the limiting process that starts at moment 0 with the
birth of a single particle. Any particle of the limiting process that
appeared at $t_{0}\in \lbrack 0,1)$ has a random lifetime uniformly
distributed on the interval $[t_{0},1)$ and, dying, generates a random
number of descendants in accordance with the generating function
\begin{equation*}
h_{\alpha }(s)=1+\left( 1+\frac{1}{\alpha }\right) \left( s-1\right) +\frac{1%
}{\alpha }\left( 1-s\right) ^{1+\alpha }.
\end{equation*}%
Sagitov's theorem generalizes the functional limit theorem due to
Fleischmann and Sigmund-Schultze \cite{FZ}, proved for reduced Galton-Watson
branching processes obeying the condition~(\ref{FinVar}).

Different questions related to the distribution of the MRCA for $k$
particles selected at random among the $Z(n)\geq k$ particles existing in
the $n$th generation were analyzed, for instance, in \cite{Athr2012}, \cite%
{Ath2012b}, \cite{Con95}, \cite{Dur78}, \cite{HJR17}, and \cite{HHP2025}.

In the present paper we consider essentially different restrictions on the $%
n $th generation population size.

\begin{theorem}
\label{T_reduced} Let condition (\ref{MainAssump}) be valid for $\alpha \in
(0,1)$ and assumption (\ref{Ratio}) hold true for a function $\varphi (n)$.
Then, for any $x\in (0,\infty )$ and any $j\geq 1$%
\begin{equation}
\lim_{n\rightarrow \infty }\mathbf{P}\left( Z(n-x\varphi (n),n)=j|\mathcal{H}%
(n)\right) =\frac{\alpha \Gamma \left( j+\alpha \right) }{j!}xM^{\ast
j}(x^{-1/\alpha }),  \label{InfVar}
\end{equation}%
where $M^{\ast j}(x)$ is the $j$th convolution of $M(x)$ with itself.
\end{theorem}

Theorem \ref{T_reduced} is a complement to the main result of \cite{LV2018},
which may be reformulated as follows (see Lemma 4 in \cite{LV2018}): if
conditions (\ref{FinVar}) and (\ref{Ratio}) are valid and the distribution
of $\xi $ has span 1, then, for any $j\geq 1$ and $x\in (0,\infty )$
\begin{equation}
\lim_{n\rightarrow \infty }\mathbf{P}\left( Z(n-x\varphi (n),n)=j\Big|%
\mathcal{H}(n)\right) =xM^{\ast j}(x^{-1}),  \label{FinVar3}
\end{equation}%
where $M(x)=1-e^{-x},x\geq 0$.

\begin{remark}
\label{R_zero} It follows from (\ref{Yaglom000}) that%
\begin{equation*}
\int_{0}^{\infty }e^{-\lambda x}dM^{\ast j}\left( x\right) =\left( 1-\frac{1%
}{\left( 1+\lambda ^{-\alpha }\right) ^{1/\alpha }}\right) ^{j}\sim \frac{1}{%
\alpha ^{j}}\frac{1}{\lambda ^{aj}},\ j=1,2,\ldots ,
\end{equation*}%
as $\lambda \rightarrow \infty $. Hence we conclude by a Tauberian theorem
(see, for instance, \cite[Ch.XIII, Sec.5, Theorems 2 and 3]{Fel}) that%
\begin{equation*}
M^{\ast j}\left( x\right) \sim \frac{1}{\alpha ^{j}\Gamma \left( 1+\alpha
j\right) }x^{\alpha j}
\end{equation*}%
as $x\rightarrow 0$. Thus,
\begin{equation}
\lim_{x\rightarrow \infty }\alpha \frac{\Gamma \left( j+\alpha \right) }{j!}%
xM^{\ast j}(x^{-1/\alpha })=\frac{\alpha \Gamma \left( j+\alpha \right) }{%
\alpha ^{j}j!\Gamma \left( 1+\alpha j\right) }\lim_{x\rightarrow \infty
}x^{1-j}=\delta _{1j},  \label{LOcal}
\end{equation}%
where $\delta _{1j}=1$ if $j=1$, and $\delta _{1j}=0$ if $j\geq 2$.
\end{remark}

Let
\begin{equation*}
b(n):=\max \left\{ 0\leq m<n:Z(m,n)=1\right\}
\end{equation*}%
be the birth moment of the MRCA of all particles existing in the population
at moment $n,$ and let $d(n):=n-b(n)$ be the distance from the moment $n$ of
observation to the birth moment of the MRCA. Since
\begin{equation*}
\mathbf{P}\left( d(n)\leq x\varphi (n)|\mathcal{H}(n)\right) =\mathbf{P}%
\left( Z(n-x\varphi (n),n)=1|\mathcal{H}(n)\right) ,
\end{equation*}%
we immediately deduce from Theorem \ref{T_reduced} the following statement.

\begin{corollary}
\label{Cor1}Under the conditions of Theorem \ref{T_reduced}, for any $x\in
(0,\infty )$%
\begin{equation*}
\lim_{n\rightarrow \infty }\mathbf{P}\left( d(n)\leq x\varphi (n)|\mathcal{H}%
(n)\right) =\alpha \Gamma \left( 1+\alpha \right) xM(x^{-1/\alpha }).
\end{equation*}
\end{corollary}

This corollary combined with (\ref{LOcal}) shows that, given $\mathcal{H}(n)$
the distance to the MRCA of all individuals of the $n$th generation is of
order $\varphi (n)$.

\section{Auxiliary results}

The proof of Theorem \ref{T_deviation} will be divided into several stages.
One of the important tools to be used in the sequel are the so-called Bell
polynomials of the second kind. To give a definition of the polynomials we
denote $\mathbb{N}:\mathbb{=}\left\{ 1,2,\ldots \right\},$ select
nonnegative integers $i_{r}\in \mathbb{N}_{0}:=\mathbb{N}\cup \left\{
0\right\} ,$ $r=1,2,\ldots ,$ and introduce, for each pair $k\leq J$ of
positive integers the set
\begin{equation*}
\Pi \left( J,k\right) :=\left\{ i_{1}+\cdots +i_{J-k+1}=k,\quad
i_{1}+2i_{2}+\cdots +(J-k+1)i_{J-k+1}=J\right\} .
\end{equation*}

The polynomials
\begin{equation}
B_{J,k}(x_{1},\ldots ,x_{J-k+1}):=\sum_{\Pi \left( J,k\right) }\frac{J!}{%
i_{1}!\cdots i_{J-k+1}!}\prod_{r=1}^{J-k+1}\left( \frac{x_{r}}{r!}\right)
^{i_{r}}  \label{DefB_Jk}
\end{equation}%
are called Bell polynomials of the second kind.

For every positive $k$, Bell polynomials of the second kind are coefficients
in the decomposition of the function
\begin{equation}
\frac{1}{k!}\left( \sum_{j=1}^{\infty }x_{j}\frac{s^{j}}{j!}\right)
^{k}=\sum_{J=k}^{\infty }B_{J,k}(x_{1},\ldots ,x_{J-k+1})\frac{s^{J}}{J!}
\label{SerB_Short}
\end{equation}%
with respect to $s$. Observe that%
\begin{equation*}
B_{J,1}(x_{1},\ldots ,x_{J})=x_{J},\text{ }B_{J,J}(x_{1})=x_{1}^{J},\text{ }
\end{equation*}%
and, for $1\leq k\leq J$%
\begin{eqnarray}
B_{J,k}(1,\ldots ,1) &=&\frac{J!}{k!}coef_{s^{J}}\left( e^{s}-1\right) ^{k}
\notag \\
&=&S(J,k)=\frac{1}{k!}\sum_{r=1}^{k}\left( -1\right) ^{k-r}\binom{k}{r}r^{J},
\label{StirlingNumbers}
\end{eqnarray}%
where $coef_{s^{J}}\left( e^{s}-1\right) ^{k}$ denotes the coefficient at $%
s^{J}$in the Taylor expansion of the function $\left( e^{s}-1\right) ^{k}$
in the vicinity of the point $s=0,$ and $S(J,k),\ k=1,2,\ldots ,J$ stand
for the Stirling numbers of the second kind. In particular, $B_{J,J}(1)=1.$

We will use Bell polynomials of the second kind to write the Fa\`{a} di
Bruno formula for the derivatives of the composition of functions $F(\cdot )$
and $G(\cdot )$ in the following form:
\begin{eqnarray*}
\frac{d^{J}}{dz^{J}}\left[ F(G(z))\right] &=&F^{\prime }(G(z))G^{(J)}(z) \\
&&+\sum_{k=2}^{J}F^{(k)}(G(z))B_{J,k}\left( G^{\prime }(z),G^{\prime \prime
}(z),\ldots ,G^{(J-k+1)}(z)\right) ,
\end{eqnarray*}%
where, as usually, the second summand is assumed to be zero if $J=1$. Hence
it follows that, for any $n,J\in \mathbb{N}$
\begin{eqnarray}
f_{n+1}^{(J)}(0) &=&\frac{d^{J}}{dz^{J}}\left[ f(f_{n}(z))\right]|_{z=0} =f^{\prime
}(f_{n}(0))f_{n}^{(J)}(0)  \notag \\
&+&\sum_{k=2}^{J}f^{(k)}(f_{n}(0))B_{J,k}\left( f_{n}^{\prime
}(0),f_{n}^{\prime \prime }(0),\ldots ,f_{n}^{(J-k+1)}(0)\right) .
\label{FaaDifference}
\end{eqnarray}

Introduce the notation%
\begin{equation*}
J_{0}:=\min \left\{ j\geq 1:\mathbf{P}\left( \xi =j\right) >0\right\} .
\end{equation*}%
Clearly,
\begin{equation}
j!\mathbf{P}\left( Z(n)=j\right) =f_{n}^{(j)}(0)=0  \label{Prob0}
\end{equation}%
for all $j=1,2,\ldots ,J_{0}-1$ and $n\in \mathbb{N}$. It follows from (\ref%
{FaaDifference}), (\ref{Prob0}) and the definition of Bell polynomials of
the second kind that
\begin{eqnarray}
J_{0}!\mathbf{P}\left( Z(n+1)=J_{0}\right) &=&f_{n+1}^{(J_{0})}(0)=f^{\prime
}(f_{n}(0))f_{n}^{(J_{0})}(0)  \notag \\
&=&\cdots =\mathbf{P}\left( \xi =J_{0}\right) \prod_{i=1}^{n}f^{\prime
}(f_{i}(0))>0  \label{ProbPositive}
\end{eqnarray}%
for each $n\in \mathbb{N}_{0}$. In fact, we can say more. First, the limit
\begin{equation}
\lim_{n\rightarrow \infty }\frac{\alpha n}{1-f_{n}(0)}\mathbf{P}\left(
Z(n)=J_{0}\right) =\mu _{J_{0}}<\infty  \label{AsymJ_0}
\end{equation}%
exists by (\ref{Stattionar}). Second, using (\ref{Stattionar2}) and
repeating almost literally the proof of Lemma 13 in \cite{NW2006} after
formula (2.41) we can conclude that $\mu _{J_{0}}>0$.

To prove Theorem \ref{T_deviation} we first analyse the asymptotic behavior
of the ratio
\begin{equation*}
q_{n}(J):=\frac{\mathbf{P}\left( Z(n)=J\right) }{\mathbf{P}\left(
Z(n)=J_{0}\right) }=\frac{J_{0}!f_{n}^{(J)}(0)}{J!f_{n}^{(J_{0})}(0)},\ J\in
\mathbb{N},
\end{equation*}%
as $n\rightarrow \infty $. It follows from (\ref{FaaDifference}) and (\ref%
{ProbPositive}) that
\begin{equation}
q_{n+1}(J)=q_{n}(J)+W_{n}(J),  \label{q_equal}
\end{equation}%
where%
\begin{equation}
W_{n}(J):=\frac{J_{0}!}{f_{n+1}^{(J_{0})}(0)}\sum_{k=2}^{J}f^{(k)}(f_{n}(0))%
\frac{B_{J,k}\left( f_{n}^{\prime }(0),f_{n}^{\prime \prime }(0),\ldots
,f_{n}^{(J-k+1)}(0)\right) }{J!}\geq 0.  \label{BellShort}
\end{equation}%
Set for brevity%
\begin{equation}
B_{J,k}:=B_{J,k}\left( 1!\mu _{1},2!\mu _{2},\ldots ,\left( J-k+1\right)
!\mu _{J-k+1}\right)  \label{Def_BellParticular}
\end{equation}%
and denote
\begin{equation}
\Delta (n):=\frac{\mathbf{P}\left( Z(n)=J_{0}\right) }{\mu _{J_{0}}}=\frac{%
f_{n}^{(J_{0})}{}(0)}{J_{0}!\mu _{J_{0}}}.  \label{Def_Delta}
\end{equation}%
Observe that
\begin{equation}
\Delta (n)\sim \frac{1-f_{n}(0)}{\alpha n}  \label{Zol}
\end{equation}%
as $n\rightarrow \infty $ according to (\ref{AsymJ_0}). \

\begin{lemma}
\label{L_w} Let condition (\ref{MainAssump}) be valid. Then there is $n_{0}$
such that
\begin{equation}
W_{n}(J)\leq \frac{2}{\mu _{J_{0}}}\sum_{k=2}^{J}f^{(k)}(f_{n}(0))\Delta
^{k-1}(n)\frac{B_{J,k}}{J!}  \label{Bell22}
\end{equation}%
for all $J\geq J_{0}$ and $n\geq n_{0}$.
\end{lemma}

\textbf{Proof}. In view of (\ref{q_equal}) and (\ref{BellShort}) the
sequence $q_{n}(j),\,n\in \mathbb{N},$ is nondecreasing for each fixed $j\in
\mathbb{N}$ and
\begin{equation*}
\lim_{n\rightarrow \infty }q_{n}(j)=\lim_{n\rightarrow \infty }\frac{%
J_{0}!f_{n}^{(j)}(0)}{j!f_{n}^{(J_{0})}(0)}=\frac{\mu _{j}}{\mu _{J_{0}}}
\end{equation*}%
according to (\ref{Stattionar}). Therefore,
\begin{equation*}
f_{n}^{(j)}(0)\leq j!\mu _{j}\frac{f_{n}^{(J_{0})}(0)}{J_{0}!\mu _{J_{0}}}%
=j!\mu _{j}\Delta (n)
\end{equation*}%
for all $j\in \mathbb{N}$ and $n\in \mathbb{N}$. It is easy to see that
\begin{eqnarray*}
&&\frac{B_{J,k}\left( f_{n}^{\prime }(0),f_{n}^{\prime \prime }(0),\ldots
,f_{n}^{(J-k+1)}(0)\right) }{J!} \\
&&\qquad \leq \frac{B_{J,k}\left( 1!\mu _{1}\Delta (n),2!\mu _{2}\Delta
(n),\ldots ,\left( J-k+1\right) !\mu _{J-k+1}\Delta (n)\right) }{J!} \\
&&\qquad \qquad =\Delta ^{k}(n)\frac{B_{J,k}}{J!}.
\end{eqnarray*}%
Hence we deduce by (\ref{BellShort}) that%
\begin{equation*}
W_{n}(J)\leq \frac{1}{\Delta \left( n+1\right) \mu _{J_{0}}}%
\sum_{k=2}^{J}f^{(k)}(f_{n}(0))\Delta ^{k}(n)\frac{B_{J,k}}{J!}.
\end{equation*}%
Since $\Delta \left( n+1\right) \sim \Delta \left( n\right) $ as $%
n\rightarrow \infty $ in view of (\ref{SurvivalProbab}), the statement of
the lemma follows.

The aim of the subsequent arguments is to estimate the difference%
\begin{eqnarray*}
0 &\leq &\sum_{J=J_{0}}^{T}\frac{\mu _{J}}{\mu _{j_{0}}}-\frac{\mathbf{P}%
\left( J_{0}\leq Z(n)\leq T\right) }{\mathbf{P}\left( Z(n)=J_{0}\right) } \\
&=&\sum_{J=J_{0}}^{T}\left( \frac{\mu _{J}}{\mu _{j_{0}}}-q_{n}(J)\right)
=\sum_{J=J_{0}}^{T}\sum_{r=n}^{\infty }\left( q_{r+1}(J)-q_{r}(J)\right) \\
&=&\sum_{J=J_{0}}^{T}\sum_{r=n}^{\infty }W_{r}(J)
\end{eqnarray*}%
from above using Lemma \ref{L_w}. It will be clear later on that to get the
desired estimates we need to evaluate the sum
\begin{equation*}
\sum_{J=k}^{T}\frac{B_{J,k}}{J!}
\end{equation*}%
as $T\geq k\rightarrow \infty $ and the derivatives $f^{(k)}(f_{n}(0)),~k%
\geq 2,$ as $n\rightarrow \infty $.

The next three sections provide such estimates.

\section{Properties of integrals containing slowly varying functions}

Recall that a function $l(x)>0,\ x\in (0,\infty ),$ is called slowly varying
at infinity if, for any $y>0$%
\begin{equation}
\lim_{x\rightarrow \infty }\frac{l(xy)}{l(x)}=1.  \label{Def_slowly}
\end{equation}

If $l(x)$ is a function slowly varying at infinity, then

1) convergence in (\ref{Def_slowly}) holds true uniformly in $y$ from any
interval $[a,b]\subset (0,\infty )$ (see \cite[Ch. 1, Sec. 1.1, Theorem 1.1]%
{Sen76});

2) for any $\delta >0$ there is $x_{0}=x_{0}(\delta )$ such that
\begin{equation}
x^{-\delta }\leq l(x)\leq x^{\delta }  \label{Slow_Approximation}
\end{equation}%
for all $x\geq x_{0}$ (see \cite[Ch.1, Sec.1.5, point $1^{\circ }$]{Sen76});

3) for any $\delta >0$ there exist $x_{0}=x_{0}(\delta )$ and $%
y_{0}=y_{0}(\delta )\geq 1$ such that, for all $x\geq x_{0}$ and $y\geq
y_{0} $
\begin{equation}
y^{-\delta }\leq \frac{l(xy)}{l(x)}\leq y^{\delta }  \label{BB2}
\end{equation}%
(see, for instance, \cite[Theorem 1.1.4 (iii)]{BB2008}).

Now we prove two lemmas concerning the asymptotic behavior of integrals
containing slowly varying functions.

\begin{lemma}
\label{L_AsymIntegralConst}Let $l(x),\ x\in (0,\infty ),$ be a function
slowly varying at infinity. Then, for any fixed $\theta >0$%
\begin{equation*}
\lim_{x\rightarrow \infty }\int_{0}^{\infty }y^{\theta -1}\frac{l(xy)}{l(x)}%
e^{-y}dy=\int_{0}^{\infty }y^{\theta -1}e^{-y}dy=\Gamma (\theta ).
\end{equation*}
\end{lemma}

\textbf{Proof}. Using (\ref{Def_slowly}), (\ref{BB2}) and convergence of the
integrals%
\begin{equation*}
\int_{0}^{\infty }y^{\theta -1\pm \delta }e^{-y}dy
\end{equation*}%
for any $\delta \in (0,\theta ),$ we get by the dominated convergence
theorem that%
\begin{eqnarray*}
\lim_{x\rightarrow \infty }\int_{0}^{\infty }y^{\theta -1}\frac{l(xy)}{l(x)}%
e^{-y}dy &=&\int_{0}^{\infty }y^{\theta -1}\lim_{x\rightarrow \infty }\frac{%
l(xy)}{l(x)}e^{-y}dy \\
&=&\int_{0}^{\infty }y^{\theta -1}e^{-y}dy=\Gamma \left( \theta \right) .
\end{eqnarray*}

Lemma \ref{L_AsymIntegralConst} is proved.

In the next lemma we consider the case $\theta \rightarrow \infty $.

\begin{lemma}
\label{L_AsumIntegral}Let $l(x),\ x\in (0,\infty ),$ be a function slowly
varying at infinity. For any $s_{0}\in (0,1)$
\begin{equation*}
\frac{\left( 1-s\right) ^{\theta }\int_{\theta }^{\infty }y^{\theta
-1}l(y)e^{-y\left( 1-s\right) }dy}{\Gamma (\theta )l\left( \frac{\theta }{1-s%
}\right) }\rightarrow 1
\end{equation*}%
as $\theta \rightarrow \infty $ uniformly in $s\in \lbrack s_{0},1)$.
\end{lemma}

\textbf{Proof}. By Stirling's asymptotic formula

\begin{equation}
\Gamma \left( \theta \right) \sim \sqrt{2\pi (\theta -1)}\left( \frac{\theta
-1}{e}\right) ^{\theta -1}\sim \sqrt{2\pi \theta }\theta ^{\theta
-1}e^{-\theta }  \label{AsymGamma}
\end{equation}%
as $\theta \rightarrow \infty $. Using the change of variables $%
y(1-s)=z\theta $ we get

\begin{eqnarray*}
\frac{1}{\Gamma \left( \theta \right) }\int_{\theta }^{\infty }y^{\theta
-1}l\left( y\right) e^{-y\left( 1-s\right) }dy &\sim &\left( \frac{e}{\theta
}\right) ^{\theta }\frac{\theta }{\sqrt{2\pi \theta }}\int_{\theta }^{\infty
}y^{\theta -1}l\left( y\right) e^{-y(1-s)}dy \\
&=&\frac{\theta e^{\theta }}{\sqrt{2\pi \theta }}\frac{1}{\left( 1-s\right)
^{\theta }}\int_{1-s}^{\infty }z^{\theta -1}l\left( \frac{\theta }{1-s}%
z\right) e^{-\theta z}dz
\end{eqnarray*}%
as $\theta \rightarrow \infty $. We fix $s_{0}\in (0,1)$ and, for $%
\varepsilon \in (0,s_{0})$ write the decomposition%
\begin{equation*}
\int_{1-s}^{\infty }z^{\theta -1}l\left( \frac{\theta }{1-s}z\right)
e^{-\theta z}dz=I_{1}(1-s,1-\varepsilon )+I_{2}(1-\varepsilon ,1+\varepsilon
)+I_{3}(1+\varepsilon ,\infty ),
\end{equation*}%
where%
\begin{equation*}
I_{k}\left( a,b\right) :=\int_{a}^{b}r(z)l\left( \frac{\theta }{1-s}z\right)
e^{\theta S(z)}dz,\quad k=1,2,3,
\end{equation*}%
and
\begin{equation*}
r(z)=\frac{1}{z},\ S(z)=\log z-z.
\end{equation*}%
Since
\begin{equation*}
\frac{l(xz)}{l(x)}\rightarrow 1
\end{equation*}%
as $x\rightarrow \infty $ uniformly in $z$ from any finite interval $%
[a,b]\subset (0,\infty )$,%
\begin{equation*}
I_{2}(1-\varepsilon ,1+\varepsilon )\sim l\left( \frac{\theta }{1-s}\right)
\int_{1-\varepsilon }^{1+\varepsilon }r(z)e^{\theta S(z)}dz
\end{equation*}%
as $\theta /\left( 1-s\right) \geq \theta /\left( 1-s_{0}\right) \rightarrow
\infty $. Clearly,%
\begin{equation*}
S^{\prime }(z)=\frac{1}{z}-1\text{ and }S^{\prime \prime }(z)=-\frac{1}{z^{2}%
}<0.
\end{equation*}%
Since $S^{\prime }(1)=0$ and \ $S^{\prime \prime }(z)<0,~z>0,$ it follows
that $z=1$ is a unique point of maximum of the function $S(z)$ on the
interval$\ z\in \lbrack 1-\varepsilon ,1+\varepsilon ]$. Therefore, we may
apply \cite[ Theorem 1.3, Sec.1, Ch.II]{Fed}, according to which%
\begin{equation*}
\int_{1-\varepsilon }^{1+\varepsilon }r(z)e^{\theta S(z)}dz\sim e^{\theta
S(1)}\sqrt{-\frac{2\pi }{\theta S^{\prime \prime }(1)}}r(1)=e^{-\theta }%
\sqrt{\frac{2\pi }{\theta }}
\end{equation*}%
as $\theta \rightarrow \infty $. Thus,
\begin{eqnarray}
\frac{\theta e^{\theta }}{\sqrt{2\pi \theta }}\frac{1}{\left( 1-s\right)
^{\theta }}I_{2}(1-\varepsilon ,1+\varepsilon ) &\sim &\frac{\theta
e^{\theta }}{\sqrt{2\pi \theta }}\frac{1}{\left( 1-s\right) ^{\theta }}%
l\left( \frac{\theta }{1-s}\right) e^{-\theta }\sqrt{\frac{2\pi }{\theta }}
\notag \\
&\sim &\frac{1}{\left( 1-s\right) ^{\theta }}l\left( \frac{\theta }{1-s}%
\right)   \label{IntermTerm}
\end{eqnarray}%
as $\theta \rightarrow \infty $ uniformly in $s\in \lbrack s_{0},1)$.

Further, in view of (\ref{BB2}), for any $\delta \in (0,1)$ there is $%
N=N\left( \delta \right) $ such that
\begin{equation*}
I_{1}(1-s,1-\varepsilon )\leq l\left( \frac{\theta }{1-s}\right)
\int_{1-s}^{1-\varepsilon }r(z)z^{-\delta }e^{\theta S(z)}dz
\end{equation*}%
and%
\begin{equation*}
I_{3}(1+\varepsilon ,\infty )\leq l\left( \frac{\theta }{1-s}\right)
\int_{1+\varepsilon }^{\infty }r(z)z^{\delta }e^{\theta S(z)}dz
\end{equation*}%
for all $\theta /\left( 1-s\right) \geq \theta /\left( 1-s_{0}\right) \geq N$%
.

Clearly,
\begin{equation*}
S(1+\varepsilon )=\log \left( 1+\varepsilon \right) -1-\varepsilon ,\quad
S^{\prime }(1+\varepsilon )=-\frac{\varepsilon }{1+\varepsilon }<0,\text{ }
\end{equation*}%
and $S^{\prime }(z)<0$ for all $z\in \lbrack 1+\varepsilon ,\infty )$.
Therefore, we may apply \cite[Theorem 1.1. Sec. 1, Ch. II]{Fed}, according
to which
\begin{eqnarray}
\int_{1+\varepsilon }^{\infty }r(z)z^{\delta }e^{\theta S(z)}dz &\sim &-%
\frac{\left( 1+\varepsilon \right) ^{\delta }r\left( 1+\varepsilon \right) }{%
\theta S^{\prime }\left( 1+\varepsilon \right) }e^{\theta S\left(
1+\varepsilon \right) }  \notag \\
&=&\frac{e^{-\theta }}{\theta }\frac{\left( 1+\varepsilon \right) ^{\delta }%
}{\varepsilon }e^{\theta \left( \log \left( 1+\varepsilon \right)
-\varepsilon \right) }.  \label{I_33}
\end{eqnarray}%
Since $\log \left( 1+\varepsilon \right) -\varepsilon <0$, it follows that%
\begin{eqnarray}
\frac{\theta e^{\theta }}{\sqrt{2\pi \theta }}\frac{1}{\left( 1-s\right)
^{\theta }}I_{3}(1+\varepsilon ,\infty ) &\leq &\frac{1}{\sqrt{2\pi \theta }}%
\frac{1}{\left( 1-s\right) ^{\theta }}l\left( \frac{\theta }{1-s}\right)
\frac{\left( 1+\varepsilon \right) ^{\delta }}{\varepsilon }e^{\theta \left(
\log \left( 1+\varepsilon \right) -\varepsilon \right) }  \notag \\
&=&o\left( \frac{1}{\left( 1-s\right) ^{\theta }}l\left( \frac{\theta }{1-s}%
\right) \right)  \label{LastTerm}
\end{eqnarray}%
as $\theta \rightarrow \infty $ uniformly in $s\in \lbrack s_{0},1)$.
Finally, consider the integral%
\begin{equation*}
\int_{1-s}^{1-\varepsilon }r(z)z^{-\delta }e^{\theta S(z)}dz.
\end{equation*}

Using the change of variables $z\rightarrow 1/w$ we get%
\begin{equation*}
\int_{1-s}^{1-\varepsilon }r(z)z^{-\delta }e^{\theta
S(z)}dz=\int_{1/(1-\varepsilon )}^{1/(1-s)}r\left( \frac{1}{w}\right)
w^{\delta -2}e^{\theta S(1/w)}dw.
\end{equation*}%
Now
\begin{equation*}
\frac{d}{dw}S\left( \frac{1}{w}\right) =\frac{d}{dw}\left( -\log w-\frac{1}{w%
}\right) =-\frac{1}{w}+\frac{1}{w^{2}}<0
\end{equation*}%
for $w\geq \left( 1-\varepsilon \right) ^{-1}$. In particular,%
\begin{equation*}
\frac{d}{dw}S\left( \frac{1}{w}\right) \left\vert _{_{w=(1-\varepsilon
)^{-1}}}\right. =-1+\varepsilon +\left( 1-\varepsilon \right)
^{2}=-\varepsilon (1-\varepsilon )<0.
\end{equation*}%
Thus, we may apply \cite[Theorem 1.1. Sec.1, Ch. II ]{Fed} once again to
conclude that
\begin{eqnarray*}
&&\int_{1/(1-\varepsilon )}^{1/(1-s)}r\left( \frac{1}{w}\right) w^{\delta
-2}e^{\theta S(1/w)}dw=\int_{1/(1-\varepsilon )}^{1/(1-s)}w^{\delta
-1}e^{-\theta (\log w+\frac{1}{w})}dw \\
&&\qquad \sim -\left( 1-\varepsilon \right) ^{\delta -1}\frac{e^{\theta
\left( \log \left( 1-\varepsilon \right) -1+\varepsilon \right) }}{\theta
\frac{d}{dw}S\left( \frac{1}{w}\right) \left\vert _{_{w=(1-\varepsilon
)^{-1}}}\right. }=\left( 1-\varepsilon \right) ^{\delta -1}\frac{e^{\theta
\left( \log \left( 1-\varepsilon \right) -1+\varepsilon \right) }}{\theta
\varepsilon \left( 1-\varepsilon \right) } \\
&&\qquad \qquad =\frac{e^{-\theta }}{\theta }\frac{1}{\left( 1-\varepsilon
\right) ^{2-\delta }}e^{\theta \left( \log \left( 1-\varepsilon \right)
+\varepsilon \right) }.
\end{eqnarray*}%
Since $\log \left( 1-\varepsilon \right) +\varepsilon <0$ for $\varepsilon
\in (0,1)$, it follows that%
\begin{eqnarray}
&&\frac{\theta e^{\theta }}{\sqrt{2\pi \theta }}\frac{1}{\left( 1-s\right)
^{\theta }}I_{1}(1-s,1-\varepsilon )  \notag \\
&&\qquad \leq \frac{\theta e^{\theta }}{\sqrt{2\pi \theta }}\frac{1}{\left(
1-s\right) ^{\theta }}l\left( \frac{\theta }{1-s}\right) \frac{e^{-\theta }}{%
\theta }\frac{1}{\left( 1-\varepsilon \right) ^{2-\delta }}e^{\theta \left(
\log \left( 1-\varepsilon \right) +\varepsilon \right) }  \notag \\
&&\qquad =\frac{1}{\sqrt{2\pi \theta }}\frac{1}{\left( 1-s\right) ^{\theta }}%
l\left( \frac{\theta }{1-s}\right) \frac{1}{\left( 1-\varepsilon \right)
^{2-\delta }}e^{\theta \left( \log \left( 1-\varepsilon \right) +\varepsilon
\right) }  \notag \\
&&\qquad =o\left( \frac{1}{\left( 1-s\right) ^{\theta }}l\left( \frac{\theta
}{1-s}\right) \right)  \label{FirstTerm}
\end{eqnarray}%
as $\theta \rightarrow \infty $ uniformly in $s\in \lbrack s_{0},1)$.

Combining (\ref{IntermTerm})--(\ref{FirstTerm}) proves the lemma.

\section{Bell polynomial of the second kind}

In this section we evaluate the sum%
\begin{equation*}
\sum_{J=k}^{T}\frac{B_{J,k}}{J!}
\end{equation*}%
as $T\geq k\geq 2$. To avoid cumbersome expressions we denote below by $%
C,C_{1},C_{2},...,$ positive constants which are not necessary the same in
different formulas.

Denote%
\begin{equation*}
H(s):=\sum_{j=1}^{\infty }\mu _{j}s^{j}.
\end{equation*}%
If condition (\ref{MainAssump}) is valid, then (see \cite[Theorem 2]{Sl1968})%
\begin{equation*}
H(s)=\frac{1}{\alpha \left( 1-s\right) ^{\alpha }L_{2}\left( 1-s\right) }%
,~0\leq s<1,
\end{equation*}%
where $L_{2}\left( z\right) \sim L\left( z\right) $ as $z\downarrow 0$.
Using (\ref{SerB_Short}) and (\ref{Def_BellParticular}) we get
\begin{equation*}
\sum_{J=k}^{\infty }\frac{B_{J,k}}{J!}s^{J}=\frac{1}{k!}\left(
\sum_{j=1}^{\infty }\mu _{j}s^{j}\right) ^{k}=\frac{1}{k!}\frac{1}{\alpha
^{k}\left( 1-s\right) ^{\alpha k}L_{2}^{k}\left( 1-s\right) }.
\end{equation*}%
Hence it follows by a Tauberian theorem (see \cite[Theorem 5, Sec. 5, Ch.
XIII]{Fel}) that, for \textbf{each fixed} $k$
\begin{equation}
\sum_{J=k}^{T}\frac{B_{J,k}}{J!}\sim \frac{1}{k!\Gamma \left( \alpha
k+1\right) }\left( \frac{T^{\alpha }}{\alpha L_{2}\left( 1/T\right) }\right)
^{k}  \label{Needed_Estimet}
\end{equation}%
as $T\rightarrow \infty $. This is not enough for us. To prove Theorem \ref%
{T_deviation} we need an estimate for the case $T\geq k\rightarrow \infty $.

To deduce the desired asymptotics in the case we interested in we introduce
the functions%
\begin{equation*}
\rho (x)=\left( \log \frac{1}{x}\right) ^{\alpha k-1},\quad 0<x<1,
\end{equation*}%
and
\begin{equation*}
g(x)=\left\{
\begin{array}{ccc}
x^{-1} & \text{if} & e^{-1}\leq x\leq 1, \\
&  &  \\
0 & \text{if} & 0\leq x<e^{-1}.%
\end{array}%
\right.
\end{equation*}

Note that%
\begin{equation}
\int_{0}^{1}\rho (x)dx=\int_{0}^{1}\left( \log \frac{1}{x}\right) ^{\alpha
k-1}dx=\int_{0}^{\infty }y^{\alpha k-1}e^{-y}dy=\Gamma \left( \alpha
k\right) ,  \label{Integ_rho}
\end{equation}%
and $g(e^{-t})=e^{t}$ if $0\leq t\leq 1$ and $g(e^{-t})=0$ if $t>1$.

The proof of the next lemma is based on the arguments used by Karamata (see
\cite{Hardy}, Ch. VII, Section 7.11, Theorems 108--110) and requires some
additional delicate estimates.

\begin{lemma}
\label{L_Bell}Let condition (\ref{MainAssump}) be valid. Then there exists $%
C>0$ such that, for any $\varepsilon >0$
\begin{equation*}
\sum_{J=k}^{T}\frac{B_{J,k}}{J!}\leq \frac{1}{k!}\frac{Ce^{\varepsilon k}}{%
\Gamma \left( \alpha k+1\right) }\left( \frac{T^{\alpha }}{\alpha
L_{2}\left( 1/T\right) }\right) ^{k}
\end{equation*}%
for all $T\geq k\geq 2$.
\end{lemma}

\textbf{Proof}. In view of (\ref{Needed_Estimet}) the statement of the lemma
is valid for all $k\leq K_{0}$ and all $T\geq k$. Thus, we may assume that $%
k $ is sufficiently large.

Introduce a non-decreasing function $\mathcal{B}_{k}(t),\ 0\leq t<\infty ,$\
by the relations
\begin{equation*}
\mathcal{B}_{k}(t):=0\text{ if }0\leq t<k,\text{ and }\mathcal{B}%
_{k}(t):=\sum_{J=k}^{\left[ t\right] }\frac{B_{J,k}}{J!},\text{ if}\ k\leq
t<\infty .
\end{equation*}%
By the definition%
\begin{eqnarray*}
\int_{0}^{\infty }e^{-\lambda t}d\mathcal{B}_{k}(t) &=&\sum_{J=k}^{\infty }%
\frac{B_{J,k}}{J!}e^{-\lambda J} \\
&=&\frac{1}{k!}H^{k}(e^{-\lambda })=\frac{1}{k!}\left( \frac{1}{\alpha
\lambda ^{\alpha }L_{3}\left( \lambda \right) }\right) ^{k},\quad \lambda >0,
\end{eqnarray*}%
where%
\begin{equation*}
L_{3}\left( \lambda \right) =\left( \frac{1-e^{-\lambda }}{\lambda }\right)
^{\alpha }L_{2}\left( 1-e^{-\lambda }\right)
\end{equation*}%
is a function slowly varying as $\lambda \downarrow 0$. Clearly, for each
fixed $m\in \mathbb{N}$ and $\varepsilon >0$ there exists $\lambda
_{0}=\lambda _{0}\left( m,\varepsilon \right) >0$ such that
\begin{eqnarray*}
\int_{0}^{\infty }e^{-\lambda t-\lambda mt}d\mathcal{B}_{k}(t) &=&\frac{1}{k!%
}\left( \frac{1}{\alpha \left( \lambda (m+1)\right) ^{\alpha }L_{3}\left(
\lambda \left( m+1\right) \right) }\right) ^{k} \\
&\leq &\frac{1}{k!}\frac{1}{\alpha ^{k}(m+1)^{\alpha k}\lambda ^{\alpha k}}%
\left( \frac{1+\varepsilon }{L_{3}\left( \lambda \right) }\right) ^{k} \\
&\leq &\frac{e^{\varepsilon k}}{k!}\left( \frac{1}{\alpha \lambda ^{\alpha
}L_{3}\left( \lambda \right) }\right) ^{k}\frac{1}{\Gamma \left( \alpha
k\right) }\int_{0}^{\infty }e^{-t}t^{\alpha k-1}e^{-mt}dt
\end{eqnarray*}%
for all $0<\lambda \leq \lambda _{0}$. Hence it follows that, for any
polynomial $R\left( x\right) ,~x\in (-\infty ,+\infty ),$ and any $%
\varepsilon >0$ there exists $\lambda _{0}=\lambda _{0}\left( R,\varepsilon
\right) >0$ such that
\begin{equation*}
\int_{0}^{\infty }e^{-\lambda t}R\left( e^{-\lambda t}\right) d\mathcal{B}%
_{k}(t)\leq \frac{e^{\varepsilon k}}{k!}\left( \frac{1}{\alpha \lambda
^{\alpha }L_{3}\left( \lambda \right) }\right) ^{k}\frac{1}{\Gamma \left(
\alpha k\right) }\int_{0}^{\infty }e^{-t}t^{\alpha k-1}R\left( e^{-t}\right)
dt
\end{equation*}%
for all $0<\lambda \leq \lambda _{0}$.

Assume that $k$ is such that%
\begin{equation*}
\Gamma ^{-3}\left( \alpha k\right) <e^{-1}\leq e^{-1}+2\Gamma ^{-3}\left(
\alpha k\right) <1.
\end{equation*}

We introduce two continuous functions on the interval $[0,1]$ by the
relations
\begin{equation*}
h_{k}^{+}(x):=\left\{
\begin{array}{ccc}
x^{-1} & \text{if} & e^{-1}\leq x\leq 1, \\
e\Gamma ^{3}\left( \alpha k\right) x+e-\Gamma ^{3}\left( \alpha k\right) &
\text{if} & e^{-1}-\Gamma ^{-3}\left( \alpha k\right) \leq x<e^{-1}, \\
0 & \text{if} & 0\leq x<e^{-1}-\Gamma ^{-3}\left( \alpha k\right) ,%
\end{array}%
\right.
\end{equation*}%
and
\begin{equation*}
h_{k}^{-}(x):=\left\{
\begin{array}{ccc}
x^{-1} & \text{if} & e^{-1}+2k\Gamma ^{-3}\left( \alpha k\right) \leq x\leq
1, \\
\frac{\Gamma ^{3}\left( \alpha k\right) \left( x-e^{-1}-\Gamma ^{3}\left(
\alpha k\right) \right) }{e^{-1}+2\Gamma ^{-3}\left( \alpha k\right) } &
\text{if} & e^{-1}+\Gamma ^{-3}\left( \alpha k\right) \leq x<e^{-1}+2\Gamma
^{-3}\left( \alpha k\right) , \\
0 & \text{if} & 0\leq x<e^{-1}+\Gamma ^{-3}\left( \alpha k\right) .%
\end{array}%
\right.
\end{equation*}%
Clearly,%
\begin{equation*}
h_{k}^{-}(x)\leq g(x)\leq h_{k}^{+}(x),~0\leq x\leq 1.
\end{equation*}%
It is not difficult to check that there is a constant $C_{1}>0$ such that%
\begin{eqnarray*}
&&\int_{0}^{1}\left( h_{k}^{+}(x)-g(x)\right) \left( \log \frac{1}{x}\right)
^{\alpha k-1}dx=\int_{e^{-1}-\Gamma ^{-3}\left( \alpha k\right)
}^{e^{-1}}h_{k}^{+}(x)\left( \log \frac{1}{x}\right) ^{\alpha k-1}dx \\
&&\qquad \qquad \qquad \leq \left( \log \frac{1}{e^{-1}-\Gamma ^{-3}\left(
\alpha k\right) }\right) ^{\alpha k-1}\int_{e^{-1}-\Gamma ^{-3}\left( \alpha
k\right) }^{e^{-1}}h_{k}^{+}(x)dx \\
&&\qquad \qquad \qquad \leq \frac{2}{\Gamma ^{3}\left( \alpha k\right) }%
\left( \log \frac{e}{1-e\Gamma ^{-3}\left( \alpha k\right) }\right) ^{\alpha
k-1} \\
&&\qquad \qquad \qquad =\frac{2}{\Gamma ^{3}\left( \alpha k\right) }\left(
1-\log \left( 1-e\Gamma ^{-3}\left( \alpha k\right) \right) \right) ^{\alpha
k-1}\leq \frac{C_{1}}{\Gamma ^{3}\left( \alpha k\right) }
\end{eqnarray*}%
and%
\begin{eqnarray*}
&&\int_{0}^{1}\left( g(x)-h_{k}^{-}(x)\right) \left( \log \frac{1}{x}\right)
^{k\alpha -1}dx=\int_{e^{-1}}^{e^{-1}+2\Gamma ^{-3}\left( \alpha k\right)
}x^{-1}\left( \log \frac{1}{x}\right) ^{k\alpha -1}dx \\
&&\qquad \qquad \qquad \qquad \leq \frac{2e}{\Gamma ^{3}\left( \alpha
k\right) }.\qquad \qquad \qquad \qquad
\end{eqnarray*}%
Since $h_{k}^{+}(x)$ is continuous, it follows from Weierstrass's theorem on
approximation of continuous functions by polynomials that there exists a
polynomial $R_{k}^{+}(x)$ such that%
\begin{equation*}
\sup_{x\in (0,1)}\left\vert R_{k}^{+}(x)-h_{k}^{+}(x)\right\vert \leq \frac{%
C_{1}}{\Gamma ^{3}\left( \alpha k\right) }.
\end{equation*}%
Taking
\begin{equation*}
P_{k}^{+}(x)=R_{k}^{+}(x)+\frac{C_{1}}{\Gamma ^{3}\left( \alpha k\right) }%
>g(x)
\end{equation*}%
we get by (\ref{Integ_rho})%
\begin{eqnarray*}
&&\int_{0}^{1}\left( P_{k}^{+}(x)-g(x)\right) \rho (x)dx\leq
\int_{0}^{1}\left( P_{k}^{+}(x)-R_{k}^{+}(x)\right) \rho (x)dx \\
&&\qquad \qquad +\int_{0}^{1}\left\vert h_{k}^{+}(x)-R_{k}^{+}(x)\right\vert
\rho (x)dx+\int_{0}^{1}\left\vert h_{k}^{+}(x)-g(x)\right\vert \rho (x)dx \\
&&\qquad \qquad \qquad \leq \frac{3C_{1}}{\Gamma ^{2}\left( \alpha k\right) }%
.
\end{eqnarray*}%
Similarly, there exists a polynomial $R_{k}^{-}(x)$ such that%
\begin{equation*}
\sup_{x\in (0,1)}\left\vert R_{k}^{-}(x)-h_{k}^{-}(x)\right\vert \leq \frac{%
C_{1}}{\Gamma ^{3}\left( \alpha k\right) }.
\end{equation*}%
Set
\begin{equation*}
P_{k}^{-}(x)=R_{k}^{-}(x)-\frac{C_{1}}{\Gamma ^{3}\left( \alpha k\right) }%
<g(x).
\end{equation*}%
Clearly,%
\begin{eqnarray*}
&&\int_{0}^{1}\left( g(x)-P_{k}^{-}(x)\right) \rho (x)dx\leq
\int_{0}^{1}\left( R_{k}^{-}(x)-P_{k}^{-}(x)\right) \rho (x)dx \\
&&\qquad \qquad \qquad +\int_{0}^{1}\left\vert
h_{k}^{-}(x)-R_{k}^{-}(x)\right\vert \rho (x)dx+\int_{0}^{1}\left\vert
h_{k}^{-}(x)-g(x)\right\vert \rho (x)dx \\
&&\qquad \qquad \qquad \qquad \leq \frac{3C_{1}}{\Gamma ^{2}\left( \alpha
k\right) }.
\end{eqnarray*}%
Now%
\begin{eqnarray*}
&&\int_{0}^{\infty }e^{-\lambda t}g\left( e^{-\lambda t}\right) d\mathcal{B}%
_{k}(t)\leq \int_{0}^{\infty }e^{-\lambda t}P_{k}^{+}\left( e^{-\lambda
t}\right) d\mathcal{B}_{k}(t) \\
&\leq &\frac{e^{\varepsilon k}}{k!}\left( \frac{1}{\alpha \lambda ^{\alpha
}L_{3}\left( \lambda \right) }\right) ^{k}\frac{1}{\Gamma \left( \alpha
k\right) }\int_{0}^{\infty }e^{-t}t^{\alpha k-1}P_{k}^{+}\left(
e^{-t}\right) dt \\
&\leq &\frac{e^{\varepsilon k}}{k!}\left( \frac{1}{\alpha \lambda ^{\alpha
}L_{3}\left( \lambda \right) }\right) ^{k}\frac{1}{\Gamma \left( \alpha
k\right) }\int_{0}^{\infty }e^{-t}t^{\alpha k-1}\left( P_{k}^{-}\left(
e^{-t}\right) +\frac{6C_{1}}{\Gamma ^{2}\left( \alpha k\right) }\right) dt \\
&=&\frac{e^{\varepsilon k}}{k!}\left( \frac{1}{\alpha \lambda ^{\alpha
}L_{3}\left( \lambda \right) }\right) ^{k}\left( \frac{1}{\Gamma \left(
\alpha k\right) }\int_{0}^{\infty }e^{-t}t^{\alpha k-1}P_{k}^{-}\left(
e^{-t}\right) dt+\frac{6C_{1}}{\Gamma ^{2}\left( \alpha k\right) }\right) \\
&\leq &\frac{e^{\varepsilon k}}{k!}\left( \frac{1}{\alpha \lambda ^{\alpha
}L_{3}\left( \lambda \right) }\right) ^{k}\left( \frac{1}{\Gamma \left(
\alpha k\right) }\int_{0}^{\infty }e^{-t}t^{\alpha k-1}g(e^{-t})dt+\frac{%
6C_{1}}{\Gamma ^{2}\left( \alpha k\right) }\right) .
\end{eqnarray*}%
Observing that
\begin{equation*}
\int_{0}^{\infty }e^{-\lambda t}g\left( e^{-\lambda t}\right) d\mathcal{B}%
_{k}(t)=\int_{0}^{1/\lambda }d\mathcal{B}_{k}(t)=\mathcal{B}_{k}(1/\lambda )
\end{equation*}%
and
\begin{eqnarray*}
\frac{1}{\Gamma \left( \alpha k\right) }\int_{0}^{\infty }e^{-t}t^{\alpha
k-1}g\left( e^{-t}\right) dt &=&\frac{1}{\Gamma \left( \alpha k\right) }%
\int_{0}^{1}t^{\alpha k-1}dt \\
&=&\frac{1}{\alpha k\Gamma \left( \alpha k\right) }=\frac{1}{\Gamma \left(
\alpha k+1\right) },
\end{eqnarray*}%
we conclude that%
\begin{equation*}
\mathcal{B}_{k}(1/\lambda )\leq \frac{e^{\varepsilon k}}{k!}\left( \frac{1}{%
\alpha \lambda ^{\alpha }L_{3}\left( \lambda \right) }\right) ^{k}\left(
\frac{1}{\Gamma \left( \alpha k+1\right) }+\frac{6C_{1}}{\Gamma ^{2}\left(
\alpha k\right) }\right) .
\end{equation*}%
Since $\Gamma \left( \alpha k+1\right) =o\left( \Gamma ^{2}\left( \alpha
k\right) \right) $ as $k\rightarrow \infty ,$ there is a constant $C>1$ such
that%
\begin{equation*}
\frac{1}{\Gamma \left( \alpha k+1\right) }+\frac{6C_{1}}{\Gamma ^{2}\left(
\alpha k\right) }\leq \frac{C}{\Gamma \left( \alpha k+1\right) }
\end{equation*}%
Using this estimate and setting $T=1/\lambda $ we conclude that
\begin{equation*}
\sum_{J=k}^{T}\frac{B_{J,k}}{J!}\leq \frac{C}{k!}\frac{e^{\varepsilon k}}{%
\Gamma \left( \alpha k+1\right) }\left( \frac{T^{\alpha }}{\alpha
L_{3}\left( 1/T\right) }\right) ^{k}.
\end{equation*}

This proves Lemma \ref{L_Bell}, since $L_{3}\left( 1/T\right) \sim
L_{2}(1/T) $ as $T\rightarrow \infty $.

\section{Asymptotics of derivatives}

In this section we evaluate the derivatives
\begin{equation*}
f^{(k)}(s),0\leq s\leq 1,\ k\geq 2,
\end{equation*}%
from above.

It is not difficult to check that if condition (\ref{MainAssump}) is valid
for $\alpha \in (0,1)$, then
\begin{equation}
p_{j}:=\mathbf{P}\left( \xi \geq j\right) =\frac{L_{1}(1/j)}{j^{1+\alpha }}%
,\quad j\geq 1,  \label{tail}
\end{equation}%
where $L_{1}(z)$ is a function slowly varying at zero and
\begin{equation}
L_{1}(z)\sim \frac{\alpha }{\Gamma \left( 1-\alpha \right) }L(z)\text{ as }%
z\downarrow 0.  \label{Slowly_ALPHA}
\end{equation}

Now we formulate the main result of this section.

\begin{lemma}
\label{L_derivative} Assume that condition (\ref{MainAssump}) is valid and $%
\alpha \in (0,1)$. Then, for any $k\geq 2$ there is $\delta >0$ such that
\begin{equation}
f^{(k)}(s)\leq k!p_{k}+2s^{-k}\frac{k\Gamma (k-1-\alpha )}{\left( 1-s\right)
^{k-1-\alpha }}L_{1}\left( \frac{1-s}{k}\right)  \label{AboveGeberal}
\end{equation}%
if $k\left( 1-s\right) \leq \delta $.
\end{lemma}

\textbf{Proof of Lemma \ref{L_derivative}}. Introduce the notation $%
x^{[k]}=x(x-1)\cdot \cdot \cdot (x-k+1),~x\geq 0,$ for falling factorials.
For $k\geq 2$ we have
\begin{eqnarray}
f^{(k)}(s) &=&\sum_{j=k}^{\infty }j^{\left[ k\right] }\mathbf{P}\left( \xi
=j\right) s^{j-k}=\sum_{j=k}^{\infty }j^{\left[ k\right] }\left(
p_{j}-p_{j+1}\right) s^{j-k}  \notag \\
&=&\sum_{j=k}^{\infty }j^{\left[ k\right] }p_{j}s^{j-k}-\sum_{j=k+1}^{\infty
}\left( j-1\right) ^{[k]}p_{j}s^{j-k-1}  \notag \\
&=&k!p_{k}+\sum_{j=k+1}^{\infty }\left( j^{\left[ k\right] }-\left(
j-1\right) ^{[k]}\right) p_{j}s^{j-k-1}  \notag \\
&&-\left( 1-s\right) \sum_{j=k+1}^{\infty }j^{\left[ k\right] }p_{j}s^{j-k-1}
\notag \\
&\leq &k!p_{k}+R(k,s),  \label{RepresentDeriv}
\end{eqnarray}%
where
\begin{equation*}
R(k,s):=k\sum_{j=k}^{\infty }j^{[k-1]}p_{j+1}s^{j-k}.
\end{equation*}%
We evaluate $R(k,s)$ from above.

Since $L_{1}(z)$ is a function slowly varying at zero, it follows that, for
any $\varepsilon _{1}>0$ there exists $x_{0}=x_{0}(\varepsilon _{1})>0$ such
that%
\begin{equation*}
\frac{L_{1}(1/\left[ x\right] )}{L_{1}(1/x)}\leq 1+\varepsilon _{1}
\end{equation*}%
for all $x\geq x_{0}.$ Using this estimate and setting $p_{x}:=p_{\left[ x%
\right] },~x>0,$ we deduce that
\begin{eqnarray*}
R(k,s) &\leq &ks^{-k-1}\int_{k}^{\infty }x^{k-1}p_{x}s^{x}dx\leq
ks^{-k-1}\int_{k}^{\max \left( x_{0}.k\right) }x^{k-1}p_{x}s^{x}dx \\
&&+(1+\varepsilon _{1})ks^{-k-1}\int_{\max \left( x_{0}.k\right) }^{\infty
}x^{k-2-\alpha }L_{1}(1/x)e^{x\log s}dx.
\end{eqnarray*}%
By the inequality $\log s\leq -(1-s),$ valid for $0<s\leq 1,$ and Lemmas \ref%
{L_AsymIntegralConst} and \ref{L_AsumIntegral} we conclude that, for any $%
\varepsilon _{2}>0$ there exists $\delta >0$ such that%
\begin{eqnarray*}
\int_{\max \left( x_{0}.k\right) }^{\infty }x^{k-2-\alpha
}L_{1}(1/x)e^{x\log s}dx &\leq &\int_{\max \left( x_{0}.k\right) }^{\infty
}x^{k-2-\alpha }L_{1}(1/x)e^{-x(1-s)}dx \\
&\leq &\frac{1+\varepsilon _{2}}{\left( 1-s\right) ^{k-1-\alpha }}%
L_{1}\left( \frac{1-s}{k}\right) \Gamma \left( k-1-\alpha \right)
\end{eqnarray*}%
as $(1-s)\max \left( x_{0}.k\right) \leq \delta $. Since the quantity%
\begin{equation*}
\int_{k}^{\max \left( x_{0}.k\right) }x^{k-1}p_{x}s^{x}dx
\end{equation*}%
is bounded and
\begin{equation*}
L_{1}\left( \frac{1-s}{k}\right) \left( 1-s\right) ^{1+\alpha -k}\rightarrow
\infty
\end{equation*}%
as $s\uparrow 1$, we conclude by the previous estimates that there exists $%
\delta >0$ such that%
\begin{equation*}
R(k,s)\leq 2s^{-k-2}\frac{k\Gamma \left( k-1-\alpha \right) }{\left(
1-s\right) ^{k-1-\alpha }}L_{1}\left( \frac{1-s}{k}\right)
\end{equation*}%
for all $(1-s)\max \left( x_{0}.k\right) \leq \delta $.

This completes the proof of Lemma \ref{L_derivative}.

\begin{lemma}
\label{L_remainder}Let condition (\ref{MainAssump}) be valid for $\alpha \in
(0,1)$. If $k\left( 1-f_{n}(0)\right) \rightarrow 0$, then
\begin{equation}
k!p_{k}=o\left( \frac{\Gamma \left( k-1-\alpha \right) }{n\left(
1-f_{n}(0)\right) ^{k-1}}\right) .  \label{O_small}
\end{equation}
\end{lemma}

\textbf{Proof}. Since $L_{1}(z)$ is a function slowly varying as $%
z\downarrow 0$, it follows that, for any $\varepsilon >0$ there is $%
z_{0}=z_{0}(\varepsilon )$ such that%
\begin{equation*}
z^{\varepsilon }\leq L_{1}(z)\leq z^{-\varepsilon }
\end{equation*}%
for all $0<z\leq z_{0}$. Hence, using (\ref{tail})\ we conclude that there
is a constant $C$ such that
\begin{equation*}
k!p_{k}\leq C\left( k-1\right) !=C\frac{\Gamma (k)n\left( 1-f_{n}(0)\right)
^{k-1}}{\Gamma \left( k-1-\alpha \right) }\times \frac{\Gamma \left(
k-1-\alpha \right) }{n\left( 1-f_{n}(0)\right) ^{k-1}}.
\end{equation*}%
for all $k\geq 2$. Thus, to prove the lemma we need to show that
\begin{equation*}
\frac{\Gamma (k)n\left( 1-f_{n}(0)\right) ^{k-1}}{\Gamma \left( k-1-\alpha
\right) }\rightarrow 0
\end{equation*}%
as $k\left( 1-f_{n}(0)\right) \rightarrow 0$. By (\ref{SimpeForm}) and (\ref%
{Slow_Approximation}) we conclude that, for any $\varepsilon \in (0,1-\alpha
)$ there is $n_{0}(\varepsilon )$ such that%
\begin{equation*}
n\leq \left( 1-f_{n}(0)\right) ^{-\alpha -\varepsilon }
\end{equation*}%
for all $n\geq n_{0}(\varepsilon )$. Thus,%
\begin{eqnarray*}
\frac{\Gamma (k)n\left( 1-f_{n}(0)\right) ^{k-1}}{\Gamma \left( k-1-\alpha
\right) } &\leq &\frac{\Gamma (k)}{\Gamma \left( k-1-\alpha \right) }\left(
1-f_{n}(0)\right) ^{k-1-\alpha -\varepsilon } \\
&=&\frac{\Gamma (k)}{\Gamma \left( k-1-\alpha \right) k^{k-1-\alpha
-\varepsilon }}\left( k(1-f_{n}(0))\right) ^{k-1-\alpha -\varepsilon }. \\
&=&o\left( \frac{\Gamma (k)}{\Gamma \left( k-1-\alpha \right) k^{k-2}}%
\right)
\end{eqnarray*}%
as $k\left( 1-f_{n}(0)\right) \rightarrow 0$. \ This proves (\ref{O_small})
for bounded $k$. Further, using (\ref{AsymGamma}) it is not difficult to
show that
\begin{equation*}
\frac{\Gamma (k)}{\Gamma \left( k-1-\alpha \right) k^{k-2}}\leq \frac{C}{%
k^{k-2}}\frac{k^{k-3/2}}{\left( k-2\right) ^{k-2-\alpha }}e^{-1-\alpha
}\rightarrow 0
\end{equation*}%
as $k\rightarrow \infty $.

The obtained estimates prove Lemma \ref{L_remainder}.

\begin{corollary}
\label{C_DerivLarge} Let condition (\ref{MainAssump}) be valid for $\alpha
\in (0,1)$. Then, for any $\gamma >0$ there exist a constant $C>0$ and $%
n_{0}=n_{0}(\gamma )$ such that
\begin{equation*}
f^{(k)}(f_{n}(0))\leq Ck^{\gamma +1}\frac{\Gamma \left( k-1-\alpha \right) }{%
n\left( 1-f_{n}(0)\right) ^{k-1}}
\end{equation*}%
for all $k\geq 2$ and all $n\geq n_{0}$ satisfying the inequality
\begin{equation}
2(k+2)(1-f_{n}(0))\leq \log 2.  \label{BasicInequal}
\end{equation}
\end{corollary}

\textbf{Proof}. By (\ref{BB2}) for any $\gamma >0$ there exist $%
k_{1}=k_{1}(\gamma )$ and $n_{1}=n_{1}(\gamma )$ such that
\begin{equation*}
\frac{L_{1}\left( \left( 1-f_{n}(0)\right) k^{-1}\right) }{L_{1}\left(
1-f_{n}(0)\right) }\leq k^{\gamma }
\end{equation*}%
for all $k\geq k_{1}$ and $n\geq n_{1}$. In view of the inequality $\log
\left( 1-x\right) \geq -2x$ valid for $x\in \lbrack 0,2^{-1}]$ we have \ $%
\log f_{n}(0)\geq -2\left( 1-f_{n}(0)\right) $ if $1-f_{n}(0)\leq 2^{-1}$
and, therefore,
\begin{equation*}
\left( f_{n}(0)\right) ^{-k-2}=e^{-(k+2)\log f_{n}(0)}\leq
e^{2(k+2)(1-f_{n}(0))}\leq 2.
\end{equation*}%
Thus,%
\begin{equation*}
\left( f_{n}(0)\right) ^{-k-2}L_{1}\left( \left( 1-f_{n}(0)\right)
k^{-1}\right) \leq 2k^{\gamma }L_{1}\left( 1-f_{n}(0)\right) .
\end{equation*}%
Now applying (\ref{AboveGeberal}) and using Lemma \ref{L_remainder} we
deduce by means of (\ref{Slowly_ALPHA}) and (\ref{SimpeForm}) that there
exist $k_{0}=k_{0}(\gamma )$ and $n_{0}=n_{0}(\gamma )$ such that
\begin{eqnarray}
&&f^{(k)}(f_{n}(0))\leq k!p_{k}+2\left( f_{n}(0)\right) ^{-k-2}\frac{k\Gamma
(k-1-\alpha )}{\left( 1-f_{n}(0)\right) ^{k-1-\alpha }}L_{1}\left( \frac{%
1-f_{n}(0)}{k}\right)  \notag \\
&&\qquad \qquad \quad \leq k!p_{k}+4k^{\gamma +1}\frac{\Gamma (k-1-\alpha )}{%
\left( 1-f_{n}(0)\right) ^{k-1-\alpha }}L_{1}\left( 1-f_{n}(0)\right)  \notag
\\
&&\qquad \qquad \quad \leq C_{1}k^{\gamma +1}\frac{\Gamma (k-1-\alpha )}{%
n\left( 1-f_{n}(0)\right) ^{k-1}}  \label{DerivLarge}
\end{eqnarray}%
for all $k\geq k_{0}$ and $n\geq n_{0}$ satisfying (\ref{BasicInequal}).

To complete the proof it remains to consider $k\in \lbrack 0,k_{0}-1]$. In
this case
\begin{equation*}
L_{1}\left( \frac{1-f_{n}(0)}{k}\right) \sim L_{1}\left( 1-f_{n}(0)\right)
\end{equation*}%
as $n\rightarrow \infty $. Thus, as before, for $k\in \lbrack 2,k_{0}-1]$
and $n\geq n_{0}$ we have%
\begin{eqnarray}
f^{(k)}(f_{n}(0)) &\leq &k!p_{k}+4\frac{k\Gamma (k-1-\alpha )}{\left(
1-f_{n}(0)\right) ^{k-1-\alpha }}L_{1}\left( \frac{1-f_{n}(0)}{k}\right)
\notag \\
&\leq &k!p_{k}+6\frac{k\Gamma (k-1-\alpha )}{\left( 1-f_{n}(0)\right)
^{k-1-\alpha }}L_{1}\left( 1-f_{n}(0)\right)  \notag \\
&\leq &k!p_{k}+C_{1}\frac{k\Gamma (k-1-\alpha )}{\left( 1-f_{n}(0)\right)
^{k-1-\alpha }}L\left( 1-f_{n}(0)\right)  \notag \\
&\leq &C_{2}k^{\gamma +1}\frac{\Gamma (k-1-\alpha )}{n\left(
1-f_{n}(0)\right) ^{k-1}}.  \label{DerivSmall}
\end{eqnarray}

Combining (\ref{DerivLarge}) and (\ref{DerivSmall}) proves Corollary \ref%
{C_DerivLarge}.

\begin{lemma}
\label{L_qEstimate} Let condition (\ref{MainAssump}) be valid for $\alpha
\in (0,1)$. Then, for any $\gamma >0$ and $\varepsilon >0$ there exist a
positive constant $C=C(\gamma ,\varepsilon )$ and $r_{0}=r_{0}(\gamma )\in
\mathbb{N}$ such that
\begin{equation*}
\sum_{k=2}^{J}f^{(k)}(f_{r}(0))\Delta ^{k-1}(r)\frac{B_{J,k}}{J!}\leq
C\sum_{k=2}^{J}\frac{k^{\gamma +1}\Gamma \left( k-1-\alpha \right) }{r^{k}}%
\left( \frac{e^{\varepsilon }}{\alpha }\right) ^{k}\frac{B_{J,k}}{J!}
\end{equation*}%
for all $J\geq 2$ and $r\geq r_{0}$ satisfying the inequality $%
2J(1-f_{r}(0))\leq \log 2$.
\end{lemma}

\textbf{Proof}. It follows from (\ref{Zol}) that, for any $\varepsilon >0$
there is $r_{1}=r_{1}(\varepsilon )$ such that
\begin{equation*}
\Delta (r)\leq (1+\varepsilon )\frac{1-f_{r}(0)}{\alpha r}\leq
e^{\varepsilon }\frac{1-f_{r}(0)}{\alpha r}
\end{equation*}%
for all $r\geq r_{1}$. Assuming that $r\geq r_{0}=\max (r_{1},n_{0})$, where
$n_{0}$ is the same as in Corollary \ref{C_DerivLarge}, we conclude that
\begin{eqnarray*}
\sum_{k=2}^{J}f^{(k)}(f_{r}(0))\Delta ^{k-1}(r)\frac{B_{J,k}}{J!} &\leq
&C\sum_{k=2}^{J}\frac{k^{\gamma +1}\Gamma \left( k-1-\alpha \right) }{%
r\left( 1-f_{r}(0\right) )^{k-1}}\left( e^{\varepsilon }\frac{1-f_{r}(0)}{%
\alpha r}\right) ^{k-1}\frac{B_{J,k}}{J!} \\
&=&C\sum_{k=2}^{J}\frac{k^{\gamma +1}\Gamma \left( k-1-\alpha \right) }{r^{k}%
}\left( \frac{e^{\varepsilon }}{\alpha }\right) ^{k-1}\frac{B_{J,k}}{J!}.
\end{eqnarray*}

Lemma \ref{L_qEstimate} is proved.

For fixed $\varepsilon ,\gamma >0$ put%
\begin{equation*}
S(J,n,\varepsilon ,\gamma )=\sum_{k=2}^{J}\frac{k^{\gamma }\Gamma \left(
k-1-\alpha \right) }{\left( n-1\right) ^{k-1}}\left( \frac{e^{\varepsilon }}{%
\alpha }\right) ^{k}\frac{B_{J,k}}{J!}.
\end{equation*}

\begin{lemma}
\label{L_muDifference} If condition (\ref{MainAssump}) holds for $\alpha \in
(0,1)$, then, for any positive $\varepsilon $ and $\gamma $ there exists a
constant $C_{1}\in (0,\infty )$ such that
\begin{equation*}
\frac{\mu _{J}}{\mu _{J_{0}}}-q_{n}(J)\leq C_{1}S(J,n,\varepsilon ,\gamma )
\end{equation*}%
for all $J\geq J_{0}$ and all  $n$ meeting the restriction $2J\left(
1-f_{n}(0)\right) \leq \log 2$.
\end{lemma}

\textbf{Proof}. By (\ref{q_equal}) and Lemmas \ref{L_w} and \ref{L_qEstimate}

\begin{eqnarray*}
q_{r+1}(J)-q_{r}(J) &\leq &\frac{2}{\mu _{J_{0}}}%
\sum_{k=2}^{J}f^{(k)}(f_{r}(0))\Delta ^{k-1}(r)\frac{B_{J,k}}{J!} \\
&\leq &C\sum_{k=2}^{J}\frac{k^{\gamma +1}\Gamma \left( k-1-\alpha \right) }{%
r^{k}}\left( \frac{e^{\varepsilon }}{\alpha }\right) ^{k}\frac{B_{J,k}}{J!}
\end{eqnarray*}%
for all $J\geq J_{0}$ and all sufficiently large $n$ such that $2J\left(
1-f_{n}(0)\right) \leq \log 2$. Since%
\begin{equation*}
\sum_{r=n}^{\infty }\frac{1}{r^{k}}\leq \frac{1}{k-1}\frac{1}{\left(
n-1\right) ^{k-1}},
\end{equation*}%
and $k\leq 2\left( k-1\right) $ for $k\geq 2$, we conclude that
\begin{equation*}
\frac{\mu _{J}}{\mu _{J_{0}}}-q_{n}(J)=\sum_{r=n}^{\infty }\left(
q_{r+1}(J)-q_{r}(J)\right) \leq 2CS(J,n,\varepsilon ,\gamma )
\end{equation*}%
for $J$ and $n$ mentioned above.

Lemma \ref{L_muDifference} is proved.

\section{\textbf{Proof of Theorem} \protect\ref{T_deviation}}

We need to evaluate the probability of the event%
\begin{equation*}
\mathcal{H}(n)=\big\{0<(1-f_{\varphi (n)}(0))Z(n)\leq 1\big\}=\big\{%
0<Z(n)\leq T\big\},
\end{equation*}%
where
\begin{equation*}
T=\left[ \frac{1}{1-f_{\varphi (n)}(0)}\right] \sim \frac{1}{1-f_{\varphi
(n)}(0)}
\end{equation*}%
as $n\rightarrow \infty $. Using \ (\ref{SurvivalProbab}) and (\ref{BB2}) it
is not difficult to verify that if $\varphi (n)n^{-1}\rightarrow 0$ as $%
n\rightarrow \infty $, then
\begin{equation*}
T\left( 1-f_{n}(0)\right) \sim \frac{1-f_{n}(0)}{1-f_{\varphi (n)}(0)}%
\rightarrow 0
\end{equation*}%
as $n\rightarrow \infty $. \ This relation and Lemma \ref{L_muDifference}
show that

\begin{equation}
\sum_{J=J_{0}}^{T}\frac{\mu _{J}}{\mu _{J_{0}}}-C_{1}S_{1}(T,n,\varepsilon
,\gamma )\leq \sum_{J=J_{0}}^{T}q_{n}(J)\leq \sum_{J=J_{0}}^{T}\frac{\mu _{J}%
}{\mu _{J_{0}}},  \label{AboveBelow}
\end{equation}%
\ where%
\begin{eqnarray}
S_{1}(T,n,\varepsilon ,\gamma ) &=&\sum_{J=J_{0}}^{T}S(J,n,\varepsilon
,\gamma )  \notag \\
&\leq &\sum_{k=2}^{T}\frac{k^{\gamma }\Gamma \left( k-1-\alpha \right) }{%
\left( n-1\right) ^{k-1}}\left( \frac{e^{\varepsilon }}{\alpha }\right)
^{k}\sum_{J=k}^{T}\frac{B_{J,k}}{J!}.  \label{S_summation}
\end{eqnarray}%
We know that $L_{2}(z)\sim L(z)$ as $z\downarrow 0$. Hence we conclude by (%
\ref{SimpeForm}) that
\begin{equation}
\frac{1}{\alpha }\frac{T^{\alpha }}{L_{2}\left( 1/T\right) }\sim \frac{1}{%
\alpha \left( 1-f_{\varphi (n)}(0)\right) ^{\alpha }L\left( 1-f_{\varphi
(n)}(0)\right) }\sim \varphi (n)  \label{MuPhi}
\end{equation}%
as $n\rightarrow \infty $. This estimate and Lemma \ref{L_Bell} lead to the
inequalities
\begin{equation*}
\sum_{J=k}^{T}\frac{B_{J,k}}{J!}\leq \frac{1}{k!}\frac{Ce^{\varepsilon k}}{%
\Gamma \left( \alpha k+1\right) }\left( \frac{T^{\alpha }}{\alpha
L_{2}\left( 1/T\right) }\right) ^{k}\leq \frac{1}{k!}\frac{Ce^{2\varepsilon
k}\varphi ^{k}(n)}{\Gamma \left( \alpha k+1\right) }
\end{equation*}%
valid for all sufficiently large $n$. Using now (\ref{S_summation}), Lemma %
\ref{L_Bell} and assumption $\varphi (n)=o(n)$ as $n\rightarrow \infty $ we
obtain for sufficiently large $n$
\begin{eqnarray*}
S_{`1}(T,n,\varepsilon ,\gamma ) &\leq &\sum_{k=2}^{T}\frac{Ce^{2\varepsilon
k}}{\left( n-1\right) ^{k-1}}\frac{k^{\gamma }\Gamma \left( k-1-\alpha
\right) }{k!\Gamma \left( \alpha k+1\right) }\varphi ^{k}(n) \\
&=&Ce^{4\varepsilon }\frac{\varphi ^{2}(n)}{n-1}\sum_{k=2}^{T}\left( \frac{%
e^{2\varepsilon }\varphi (n)}{n-1}\right) ^{k-2}\frac{k^{\gamma }\Gamma
\left( k-1-\alpha \right) }{k!\Gamma \left( \alpha k+1\right) } \\
&\leq &Ce^{4\varepsilon }\frac{\varphi ^{2}(n)}{n-1}\sum_{k=2}^{T}\frac{%
k^{\gamma }\Gamma \left( k-1-\alpha \right) }{k!\Gamma \left( \alpha
k+1\right) }.
\end{eqnarray*}%
It is easy to see by Stirling's asymptotic formula that the series%
\begin{equation*}
\sum_{k=2}^{\infty }\frac{k^{\gamma }\Gamma \left( k-1-\alpha \right) }{%
k!\Gamma \left( \alpha k+1\right) }
\end{equation*}%
is convergent. Since $\varphi (n)=o(n)$ as $n\rightarrow \infty $, it
follows that $S_{1}(T,n,\varepsilon ,\gamma )=o\left( \varphi (n)\right) $
as $n\rightarrow \infty $. This estimate combined with (\ref{Sl2}), (\ref%
{AboveBelow}) and (\ref{MuPhi}) shows that%
\begin{equation*}
\sum_{J=J_{0}}^{T}q_{n}(J)\sim \sum_{J=J_{0}}^{T}\frac{\mu _{J}}{\mu _{J_{0}}%
}\rightarrow \infty
\end{equation*}%
as $n\rightarrow \infty $. As a result,
\begin{eqnarray*}
\mathbf{P}\left( 0<Z(n)\leq T\right)  &=&\mathbf{P}\left( Z(n)=J_{0}\right)
\sum_{J=1}^{T}\frac{\mathbf{P}\left( Z(n)=J\right) }{\mathbf{P}\left(
Z(n)=J_{0}\right) } \\
&\sim &\frac{\mathbf{P}\left( Z(n)=J_{0}\right) }{\mu _{J_{0}}}%
\sum_{J=1}^{T}\mu _{J}\sim \Delta \left( n\right) \frac{1}{\alpha \Gamma
\left( 1+\alpha \right) }\frac{T^{\alpha }}{L\left( 1/T\right) } \\
&\sim &\Delta \left( n\right) \frac{\varphi (n)}{\Gamma \left( 1+\alpha
\right) }\sim \frac{1-f_{n}(0)}{\alpha n}\frac{\varphi (n)}{\Gamma \left(
1+\alpha \right) }.
\end{eqnarray*}

Theorem \ref{T_deviation} is proved.

\begin{remark}
\label{R_alpha1} As it was mentioned in section 1, we do not consider the case $%
\alpha =1$ and $\sigma ^{2}=\infty $. Simple transformation show that in
this case
\begin{equation*}
\frac{f(s)-s}{\left( 1-s\right) ^{2}}=\sum_{r=0}^{\infty }\left(
\sum_{j=r+1}^{\infty }\mathbf{P}\left( \xi \geq j\right) \right)
s^{r}=L(1-s),
\end{equation*}%
where $L(z)\rightarrow \infty $ as $z\downarrow 0$. Hence it follows that,
as $n\rightarrow \infty $%
\begin{equation*}
\sum_{r=0}^{n}\left( \sum_{j=r+1}^{\infty }\mathbf{P}\left( \xi \geq
j\right) \right) \sim L\left( 1/n\right) .
\end{equation*}%
Unfortunately, we cannot find from this relation the asymptotic behavior of
the sum
\begin{equation*}
\sum_{j=n}^{\infty }\mathbf{P}\left( \xi \geq j\right)
\end{equation*}%
and the probabilities $p_{n}=\mathbf{P}\left( \xi \geq n\right) $ as $%
n\rightarrow \infty $. However, if we assume that%
\begin{equation}
\mathbf{P}\left( \xi \geq j\right) =\frac{\bar{L}(1/j)}{j^{2}},
\label{Alpha_1}
\end{equation}%
where $\bar{L}(z)$ is a function slowly varying at zero, then the
representation
\begin{eqnarray*}
\sum_{r=0}^{n}\left( \sum_{j=r+1}^{\infty }\mathbf{P}\left( \xi \geq
j\right) \right)  &=&\sum_{r=0}^{n}(1+\varepsilon _{1}(r))\frac{\bar{L}(1/r)%
}{r+1} \\
&=&\sum_{r=0}^{n}(1+\varepsilon _{2}(r))rp_{r+1}\sim L\left( 1/n\right) ,
\end{eqnarray*}%
where $\left\vert \varepsilon _{1}(r)\right\vert +\left\vert \varepsilon
_{2}(r)\right\vert \rightarrow 0$ as $r\rightarrow \infty $ implies (see
\cite{Fel}, Ch. VIII, Sec.9, formula (9.5) for $p=0,\gamma =-1$) that%
\begin{equation}
\bar{L}\left( z\right) =o\left( L(z)\right) =o\left( L_{1}(z)\right)
\label{Negl_L}
\end{equation}
as $z\downarrow 0$. Using this estimate and (\ref{RepresentDeriv}) we
conclude by \cite[Ch. XIII, Sec. 5, Theorem 5]{Fel} that
\begin{equation*}
f^{(2)}(s)\leq 2!p_{2}+2\sum_{j=3}^{\infty }jp_{j+1}s^{j-2}\leq CL\left(
1-s\right) .
\end{equation*}%
Now relation (\ref{SimpeForm})\ leads to the estimate%
\begin{eqnarray*}
f^{(2)}(f_{n}(0)) &\leq &CL\left( 1-f_{n}(0)\right) =C\frac{L\left(
1-f_{n}(0)\right) \left( 1-f_{n}(0)\right) }{\left( 1-f_{n}(0)\right) } \\
&\leq &\frac{C_{1}}{n\left( 1-f_{n}(0)\right) }.
\end{eqnarray*}%
Moreover, if (\ref{Alpha_1}) is valid, then the estimates established in
Lemmas \ref{L_derivative}, \ref{L_remainder} and Corollary \ref{C_DerivLarge}
are valid for $\alpha =1$ and $k\geq 3,$ where $L_{1}(\cdot )$ is replaced
by $\bar{L}\left( \cdot \right) $. Using this fact and estimate (\ref{Negl_L}%
) one can show that the conclusion of Theorem~\ref{T_deviation} remains
valid, if assumption (\ref{Alpha_1}) holds and $\sigma ^{2}=\infty $.
\end{remark}

\section{Limit theorem for reduced processes}

In this section we will prove Theorem \ref{T_reduced}.

Let $m=n-x\varphi (n)$, where $x\in \left( 0,\infty \right) $ and $\varphi
(n)\rightarrow \infty ,\varphi (n)=o(n)$ as $n\rightarrow \infty $.

We start by checking the following statement.

\begin{lemma}
\label{L_derivativeLarge} Let condition (\ref{MainAssump}) be valid. Then,
for any fixed $J=1,2,\ldots $
\begin{equation*}
f_{m}^{(J)}\left( f_{x\varphi (n)}(0)\right) \sim \left( -1\right) ^{J}\frac{%
x\varphi (n)\Delta (n)}{\left( \log f_{x\varphi (n)}(0)\right) ^{J}}\frac{%
\Gamma \left( J+\alpha \right) }{\Gamma \left( \alpha \right) }\text{ as }%
n\rightarrow \infty \text{.}
\end{equation*}
\end{lemma}

\textbf{Proof}. It is known (see, for instance, \cite[Lemmas 2 and 4]{Sl1968}%
) that under condition (\ref{MainAssump})
\begin{eqnarray*}
&&\lim_{n\rightarrow \infty }\left( n\alpha \right) ^{1+1/\alpha }L\left(
1-f_{n}(0)\right) ^{1/\alpha }\left( f_{n+1}(0)-f_{n}(0)\right) \\
&&\qquad =\lim_{n\rightarrow \infty }n\alpha \left( \frac{f_{n+1}(0)-f_{n}(0)%
}{1-f_{n}(0)}\right) =\lim_{n\rightarrow \infty }\frac{f_{n+1}(0)-f_{n}(0)}{%
\Delta (n)}=1,
\end{eqnarray*}%
where we have used (\ref{Def_Delta}) and (\ref{Zol}).

We consider the function%
\begin{equation*}
f_{m}(f_{x\varphi (n)}^{\lambda }(0))=f_{m}(e^{\lambda \log f_{x\varphi
(n)}(0)}),\ \lambda >0,
\end{equation*}%
and, for each $\lambda >0$ find a parameter $r=r(\lambda )$ such that%
\begin{equation*}
1-f_{r+1}(0)<1-f_{x\varphi (n)}^{\lambda }(0)\leq 1-f_{r}(0).
\end{equation*}%
In view of (\ref{SurvivalProbab})
\begin{eqnarray*}
1-f_{x\varphi (n)}^{\lambda }(0) &\sim &\lambda \left( 1-f_{x\varphi
(n)}(0)\right) \sim \frac{\lambda L^{\ast }(x\varphi (n))}{\left( x\varphi
(n)\right) ^{1/\alpha }} \\
&\sim &\frac{L^{\ast }(x\varphi (n)/\lambda ^{a})}{\left( x\varphi
(n)/\lambda ^{\alpha }\right) ^{1/\alpha }}\sim \frac{L^{\ast }(r)}{%
r^{1/\alpha }}
\end{eqnarray*}%
as $n\rightarrow \infty $. Hence we get%
\begin{equation*}
r\sim \frac{x\varphi (n)}{\lambda ^{\alpha }}=o(n)\text{ as }n\rightarrow
\infty \text{.}
\end{equation*}%
Since $m=n-x\varphi (n)\sim n,$ we have%
\begin{eqnarray*}
&&\lim_{n\rightarrow \infty }\frac{1}{x\varphi (n)}\left[ \frac{%
f_{m}(f_{r}(0))-f_{m}(0)}{\Delta \left( n\right) }\right] =\lim_{n%
\rightarrow \infty }\frac{1}{x\varphi (n)}\sum_{k=1}^{r-1}\frac{%
f_{m}(f_{k+1}(0))-f_{m}(f_{k}(0))}{\Delta \left( n\right) } \\
&&\qquad \qquad \qquad \qquad \qquad =\lim_{n\rightarrow \infty }\frac{1}{%
x\varphi (n)}\sum_{k=1}^{r-1}\frac{\Delta \left( m+k+1\right) }{\Delta
\left( n\right) }\left[ \frac{f_{m+k+1}(0)-f_{m+k}(0)}{\Delta \left(
m+k+1\right) }\right] \\
&&\qquad \qquad \qquad \qquad \qquad \qquad =\lim_{n\rightarrow \infty }%
\frac{1}{x\varphi (n)}\sum_{k=1}^{r-1}1=\frac{1}{\lambda ^{\alpha }}.
\end{eqnarray*}%
Thus,%
\begin{equation}
\lim_{n\rightarrow \infty }\frac{f_{m}(e^{\lambda \log f_{x\varphi
(n)}(0)})-f_{m}(0)}{x\varphi (n)\Delta \left( n\right) }=\frac{1}{\lambda
^{\alpha }},\quad \lambda >0.  \label{Complex}
\end{equation}%
Since the prelimiting functions at the left-hand side of (\ref{Complex}) are
analytical in the domain $Re \lambda >0$, it follows that%
\begin{equation}
\lim_{n\rightarrow \infty }\frac{1}{x\varphi (n)\Delta \left( n\right) }%
\frac{d^{J}}{d^{J}\lambda }\left[ f_{m}(e^{\lambda \log f_{x\varphi (n)}(0)})%
\right] =\left( -1\right) ^{J}\frac{\Gamma (\alpha +J)}{\Gamma (\alpha )}%
\frac{1}{\lambda ^{\alpha +J}},\ \lambda >0,  \label{Derivative2}
\end{equation}%
for each $J\geq 1$. By (\ref{FaaDifference}), (\ref{DefB_Jk}) and (\ref%
{StirlingNumbers}) we have
\begin{eqnarray*}
&&\frac{d^{J}}{d^{J}\lambda }\left[ f_{m}(e^{\lambda \log f_{x\varphi
(n)}(0)})\right] =\sum_{k=1}^{J}f_{m}^{(k)}(f_{x\varphi (n)}^{\lambda }(0))
\\
&&\times B_{J,k}\left( f_{x\varphi (n)}^{\lambda }(0)\log f_{x\varphi
(n)}(0),f_{x\varphi (n)}^{\lambda }(0)\log ^{2}f_{x\varphi (n)}(0),\ldots
,f_{x\varphi (n)}^{\lambda }(0)\left( \log f_{x\varphi (n)}(0)\right)
^{J-k+1}\right) \\
&&\qquad =\left( \log f_{x\varphi (n)}(0)\right)
^{J}\sum_{k=1}^{J}f_{m}^{(k)}(f_{x\varphi (n)}^{\lambda }(0))\left(
f_{x\varphi (n)}^{\lambda }(0)\right) ^{k}B_{J,k}\left( 1,1,\ldots ,1\right)
\\
&&\qquad =\left( \log f_{x\varphi (n)}(0)\right)
^{J}\sum_{k=1}^{J}f_{m}^{(k)}(f_{x\varphi (n)}^{\lambda }(0))\left(
f_{x\varphi (n)}^{\lambda }(0)\right) ^{k}S(J,k).
\end{eqnarray*}%
It follows from this representation and (\ref{Derivative2}) that%
\begin{equation*}
\frac{1}{x\varphi (n)\Delta \left( n\right) }\frac{d}{d\lambda }\left[
f_{m}(e^{\lambda \log f_{x\varphi (n)}(0)})\right] \left\vert _{\lambda
=1}\right. \sim \frac{\log f_{x\varphi (n)}(0)}{x\varphi (n)\Delta \left(
n\right) }f_{m}^{\prime }(f_{x\varphi (n)}(0))\sim \left( -1\right) \alpha
\end{equation*}%
as $n\rightarrow \infty $, giving
\begin{equation*}
f_{m}^{\prime }(f_{x\varphi (n)}(0))\sim \left( -1\right) \alpha \frac{%
x\varphi (n)\Delta \left( n\right) }{\log f_{x\varphi (n)}(0)}.
\end{equation*}%
Now we prove by induction that, for any $k\geq 1$
\begin{equation*}
\frac{\left( \log f_{x\varphi (n)}(0)\right) ^{k}}{x\varphi (n)\Delta \left(
n\right) }f_{m}^{(k)}\left( f_{x\varphi (n)}(0)\right) \sim \left( -1\right)
^{k}\frac{\Gamma (\alpha +k)}{\Gamma (\alpha )}
\end{equation*}%
as $n\rightarrow \infty $. This is true for $k=1$ and if this is true for
all $k<J$ then, in view of (\ref{Derivative2}) and convergence of $\log
f_{x\varphi (n)}(0)$ to 0 as $n\rightarrow \infty ,$ we conclude by
induction hypothesis that
\begin{eqnarray*}
&&\frac{1}{x\varphi (n)\Delta \left( n\right) }\frac{d^{J}}{d^{J}\lambda }%
\left[ f_{m}(e^{\lambda \log f_{x\varphi (n)}(0)})\right] \left\vert
_{\lambda =1}\right. =\frac{\left( \log f_{x\varphi (n)}(0)\right)
^{J}f_{m}^{(J)}(f_{x\varphi (n)}(0))}{x\varphi (n)\Delta \left( n\right) } \\
&&\qquad \qquad \qquad \qquad \qquad \quad \quad +\sum_{k=1}^{J-1}\frac{%
\left( \log f_{x\varphi (n)}(0)\right) ^{J}f_{m}^{(k)}(f_{x\varphi (n)}(0))}{%
x\varphi (n)\Delta \left( n\right) }S(J,k) \\
&&\qquad \qquad \qquad \qquad =\frac{\left( \log f_{x\varphi (n)}(0)\right)
^{J}f_{m}^{(J)}(f_{x\varphi (n)}(0))}{x\varphi (n)\Delta \left( n\right) }%
+O\left( \log f_{x\varphi (n)}(0)\right) \\
&&\qquad \qquad \qquad \qquad \sim \left( -1\right) ^{J}\frac{\Gamma (\alpha
+J)}{\Gamma (\alpha )}
\end{eqnarray*}%
as $n\rightarrow \infty $.

Lemma \ref{L_derivativeLarge} is proved.

\textbf{Proof of Theorem \ref{T_reduced}}. For any $j\geq 1$ and $%
m=n-x\varphi (n)$ we have%
\begin{eqnarray*}
\mathbf{P}\left( Z(m,n\right) =j) &=&\sum_{k=j}^{\infty }\mathbf{P}\left(
Z(m)=k;Z(m,n\right) =j) \\
&=&\sum_{k=j}^{\infty }\mathbf{P}\left( Z(m)=k\right) \binom{k}{j}%
f_{n-m}^{k-j}(0)\left( 1-f_{n-m}(0)\right) ^{j} \\
&=&\frac{\left( 1-f_{n-m}(0)\right) ^{j}}{j!}f_{m}^{(j)}(f_{n-m}(0)).
\end{eqnarray*}%
This representation and Lemma \ref{L_derivativeLarge} show that, as $%
n\rightarrow \infty $
\begin{eqnarray*}
\mathbf{P}\left( Z(n-x\varphi (n),n\right)  &=&j)=\frac{\left( 1-f_{x\varphi
(n)}(0)\right) ^{j}}{j!}f_{n-x\varphi (n)}^{(j)}\left( f_{x\varphi
(n)}(0)\right)  \\
&\sim &\left( -1\right) ^{j}\frac{\left( 1-f_{x\varphi (n)}(0)\right) ^{j}}{%
j!}\frac{x\varphi (n)\Delta \left( n\right) }{\left( \log f_{x\varphi
(n)}(0)\right) ^{j}}\frac{\Gamma (\alpha +j)}{\Gamma (\alpha )} \\
&\sim &\frac{\Gamma (\alpha +j)}{\Gamma (\alpha )j!}x\varphi (n)\Delta
\left( n\right) .
\end{eqnarray*}%
Denote $Z_{1}^{\ast }(l),\ldots ,Z_{j}^{\ast }(l)$ i.i.d. random variables
distributed as $\left\{ Z(l)|Z(l)>0\right\} ,$ and let $\eta _{1},\ldots
,\eta _{j}$ be i.i.d. random variables having distribution $M(x)$ specified
by (\ref{Yaglom000}). Using (\ref{Def_H}), (\ref{Yaglom000}) and (\ref%
{SurvivalProbab}) it is not difficult to check that
\begin{eqnarray*}
&&\lim_{n\rightarrow \infty }\mathbf{P}\left( \mathcal{H}(n)|Z(m,n)=j\right)
=\lim_{n\rightarrow \infty }\mathbf{P}\left( \sum_{r=1}^{j}Z_{r}^{\ast
}(x\varphi (n))\leq \frac{1}{1-f_{\varphi (n)}(0)}\right)  \\
&&\qquad \qquad =\lim_{n\rightarrow \infty }\mathbf{P}\left(
\sum_{r=1}^{j}\left( 1-f_{x\varphi (n)}(0)\right) Z_{r}^{\ast }(x\varphi
(n))\leq \frac{1-f_{x\varphi (n)}(0)}{1-f_{\varphi (n)}(0)}\right)  \\
&&\qquad \qquad \qquad =\mathbf{P}\left( \eta _{1}+\cdots +\eta _{j}\leq
x^{-1/\alpha }\right) =M^{\ast j}(x^{-1/\alpha }).
\end{eqnarray*}%
Combining this result with Lemma \ref{L_derivativeLarge} and Theorem \ref%
{T_deviation} we see that%
\begin{eqnarray*}
&&\mathbf{P}\left( Z(m,n)=j|\mathcal{H}(n)\right) =\frac{\mathbf{P}\left(
Z(m,n\right) =j)\mathbf{P}\left( \mathcal{H}(n)|Z(m,n)=j\right) }{\mathbf{P}%
\left( \mathcal{H}(n)\right) } \\
&&\qquad \qquad \qquad \quad \sim \frac{\Gamma (\alpha +j)}{\Gamma (\alpha
)j!}x\varphi (n)\Delta \left( n\right) \frac{\Gamma \left( 1+\alpha \right)
}{\Delta \left( n\right) \varphi (n)}M^{\ast j}(x^{-1/\alpha }) \\
&&\qquad \qquad \qquad \qquad \quad \qquad =\frac{\alpha \Gamma (\alpha +j)}{%
j!}xM^{\ast j}(x^{-1/\alpha }).
\end{eqnarray*}%
Theorem \ref{T_reduced} is proved.

\begin{remark}
The conclusion of Theorem \ref{T_reduced} remains valid, if
assumption (\ref{Alpha_1}) holds and $\sigma ^{2}=\infty $ (see Remark \ref%
{R_alpha1}).
\end{remark}

\hfill \rule{2mm}{3mm}

\end{document}